\documentclass[11pt,twoside]{article}
\usepackage{latexsym}
\usepackage{bbold, epsfig}
\usepackage{amssymb,amsbsy,amsmath,amsfonts,amssymb,amscd}
\usepackage[T1]{fontenc}
\usepackage{indentfirst}
\usepackage{fancyhdr, fullpage, enumerate,array}

\usepackage{ncclatex}
\AfterTheoremHeaderChar{.}

\newtheorem{theo}{Theorem}

\newtheorem{lem}{Lemma}
\newtheorem{cor}{Corollary}

\newtheorem*{defi*}{Definition}
\newtheorem{rem}{Remark}

\newcommand{\D}{\mbox \DH}
\newcommand{\ds}{\displaystyle}
\renewcommand{\div}{{\rm div \,}}
\newcommand{\grad}{{\rm grad \,}}
\newcommand{\R}{\mathbb {R}}
\newcommand{\N}{\mathbb {N}}

\newcommand{\divv}{\operatorname{div}}
\newcommand{\inj}{\hookrightarrow}
\newcommand{\ep}{\varepsilon}

\newcommand{\nd}{\noindent}

\newcommand{\cqfd}{\begin{flushright}\vspace*{-3mm}$\Box $\vspace{-2mm}\end{flushright}}

\newcommand{\on}{\mbox{ on }}

\newcommand{\eps}{\varepsilon}

\newcommand{\saut}{\vspace*{1mm}\\ \nd }

\DeclareMathOperator{\esssup}{ess\,sup}

\begin{document}
\title{Global existence for a class of reaction-diffusion systems with mass action kinetics and concentration-dependent diffusivities}
\author{Dieter Bothe,\; Guillaume Rolland\\
}
\maketitle
%

\begin{abstract}
In this work we study the existence of classical solutions for a class of reaction-diffusion systems with quadratic growth naturally arising in mass action chemistry when studying networks of reactions of the type $A_i+A_j \rightleftharpoons A_k$ with Fickian diffusion, where the diffusion coefficients might depend on time, space and on all the concentrations $c_i$ of the chemical species. In the case of one single reaction, we prove global existence for space dimensions $N\leq 5$. In the more restrictive case of diffusion coefficients of the type $d_i(c_i)$, we use an $L^2$-approach to prove global existence for $N\leq 9$.
In the general case of networks of such reactions we extend the previous method to get global solutions  for general diffusivities if $N\leq 3$ and for diffusion of type $d_i(c_i)$ if $N\leq 5$.
In the latter quasi-linear case of $d_i(c_i)$ and for space dimensions $N=2$ and $N=3$, global existence holds for more than quadratic reactions. We can actually allow for more general rate functions including fractional power terms, important in applications. We obtain global existence under
appropriate growth restrictions with an explicit dependence on the space dimension $N$.\\
 \end{abstract}



\section{Introduction}
Chemical reaction-diffusion systems (RD-systems for short) consist of mass balances, often given in terms of molar mass densities
$c_i$ of certain chemical species $A_i$, where $i=1,\ldots ,P$ in case of $P$ involved chemical components.
This leads to PDE-systems of the form
\begin{equation}
\label{mass-balance}
\partial_t c_i + \div J_i = f_i \qquad (i=1, \ldots ,P),
\end{equation}
where $J_i$ is the (molar) mass flux of species $A_i$ and the source term $r_i$ models the rate of change of $A_i$ due to chemical reactions.

While transport of $A_i$ is usually mediated by several parallel mechanisms like convection, diffusion or migration,
the fluxes in \eqref{mass-balance} are commonly considered to be of diffusive type in case of RD-systems.
These diffusive fluxes are most often modeled by
the classical Fick's law, i.e.\ constitutive relations of the type
\begin{equation}
\label{Ficks-law}
J_i = - d_i \, \nabla c_i  \qquad (i=1, \ldots ,P)
\end{equation}
are employed for this purpose, where the diffusivities $d_i$ are nonnegative due to the second law of thermodynamics
\cite{deGrootMazur-book}. In \eqref{Ficks-law}, the $d_i$ will be (complicated) functions of the system's thermodynamic state variables,
in particular the diffusivities depend significantly on the mixture composition, i.e.\ on the concentration vector $c:=(c_1, \ldots ,c_P)$.
A flux of Fickian type \eqref{Ficks-law} can either model so-called molecular diffusion caused by the random thermal motion of all molecules, or an effective diffusive flux due to other stochastic particle motions such as random convective motions of fluid parcels in a turbulent velocity field. In the later case one also speaks of dispersive mixing or dispersion instead of diffusion; cf.\ \cite{BaldygaBourne}.

We consider systems of type \eqref{mass-balance} in bounded domains $\Omega \subset \R^N$ with sufficiently smooth boundary $\partial \Omega$
under the homogeneous Neumann boundary condition
\begin{equation}
\label{Neumann-BC}
J_i \cdot \nu = 0 \quad \mbox{ on $\partial \Omega\qquad$ } (i=1, \ldots ,P),
\end{equation}
where $\nu$ denotes the outward unit normal to $\Omega$.
We also impose the initial conditions
\begin{equation}
\label{ICs}
c_i (0,\cdot)  = c_{0,i} \quad \mbox{ on $\Omega\qquad$ } (i=1, \ldots ,P),
\end{equation}
where the initial concentrations $c_{0,i}$ are nonnegative and sufficiently regular, at least $L^\infty$.

One main emphasis of the present paper lies on the investigation of RD-systems from physico-chemical backgrounds.
Typical applications come from Chemical Reaction Engineering, say reactions in liquid systems under isobaric conditions
(such that no convective flow occurs) or diffusion of reactive species into solids.
There is a large amount of measurement data from such applications, showing the dependence of the Fickian diffusivities on the
concentrations; see in particular \cite{cussler}.
Instead of going into further details on measured dependencies, we prefer to include a brief theoretical explanation.
For such systems, the Maxwell-Stefan equations provide a more fundamental and thermodynamically consistent approach to model
diffusive multicomponent transport; cf.\ \cite{Giovan}, \cite{TaylorKrishna}.
The  Maxwell-Stefan equations form a reduced set of partial momentum balances for the involved
constituents, relying on a scale-separation argument which is a very accurate approximation for diffusion velocities far below
the speed of sound \cite{BotheDreyer}.
To avoid cases with additional migrative transport, we also assume that the species are uncharged, which rules
out certain cases with ionic species especially appearing in aqueous solutions.
Furthermore, we assume isothermal conditions to avoid thermo-diffusion processes and, more important,
severe complications due to the usually significant temperature dependence of chemical reactions.
Finally, we assume that no convective transport occurs in the mixture. In case of a fluid system this corresponds to isobaric conditions,
since any pressure gradient will cause the mixture to flow.
In the resulting isobaric and isothermal case without species-dependent body forces, the Maxwell-Stefan equations read
\begin{equation}
\label{MS-system}
 -\, \sum_{j\neq i} \frac{ x_j {\, J}_i - x_i {\, J}_j}{c_{\rm tot}\, \D_{ij}} =  \frac{x_i}{R\, T}\, \nabla \mu_i
\quad \mbox{ for } i=1,\ldots ,P.
\end{equation}
Here $c_{\rm tot}:=\sum_i c_i$ is the total concentration, $x_i:=c_i/c_{\rm tot}$ are the molar fractions, $R$ is the
universal gas constant, $T$ the absolute temperature and $\mu_i$ the chemical potential of species $A_i$. Moreover, the
$\D_{ij}$ are the so-called Maxwell-Stefan diffusivities which are symmetric, where the latter is either seen as a consequence
of Onsager's reciprocal relations, or can be deduced under the assumption of binary interactions; cf.\ \cite{BotheDreyer}.
Like the Fickian diffusivities, the $\D_{ij}$ are not constant but depend on the thermodynamic state variables
 - especially, $\D_{ij}=\D_{ij} (c)$.
The set of equations \eqref{MS-system} is complemented by the constraint
\begin{equation}
\label{flux-sum}
\sum_{i=1}^P {\, J}_i =0,
\end{equation}
expressing the fact that diffusive fluxes are taken relative to a common, molar
averaged mixture velocity, where the latter is assumed to be zero throughout this paper.
For a more complete model with consistent coupling to the {\it barycentric} momentum balance, see \cite{BotheDreyer}.

The system of equations \eqref{MS-system} and \eqref{flux-sum} can be inverted to obtain the diffusive fluxes $J_i$; see \cite{Giovan}, \cite{BoMS}.
The resulting fluxes account both for direct cross-effects due to friction between the components as expressed by the left-hand side
in \eqref{MS-system}, and for non-idealities due to complex material behavior which enters via the chemical potentials on the
right-hand side of \eqref{MS-system}.
In the general case of a multicomponent system with diffusive fluxes modeled by \eqref{MS-system} and \eqref{flux-sum}, a
fully coupled RD-system with fluxes of type
\begin{equation}
\label{Ficks-law-gen}
J_i = - \sum_{j=1}^P d_{ij} \, \nabla c_j  \qquad (i=1, \ldots ,P)
\end{equation}
results, where the non-diagonal diffusion matrix $[d_{ij}]$ depends on the composition $c$.
Without chemical reactions, the pure diffusion system \eqref{mass-balance}, \eqref{Neumann-BC} -- \eqref{flux-sum}
is locally in time wellposed for sufficiently regular initial data as shown in \cite{BoMS}. But for the chemically reactive case, only first results on existence of global weak solutions (without uniqueness) are available in \cite{Juengel} for rather restricted chemistry.

The present paper investigates the complications due to non-constant diffusivities, but possible diffusive cross-effects are ignored.
To motivate these particular class of RD-systems with concentration-dependent diffusivities but without cross-diffusion, let us
briefly discuss two important special cases in which the Maxwell-Stefan equations can be explicitly inverted.
For a \emph{binary mixture}, i.e.\ a mixture with two components, it follows from
$x_1 + x_2=1$ and ${ J}_1+{ J}_2=0$ that
\begin{equation}
\label{binary2}
{ J}_1 \, ( =- { J}_2 )\, = -\frac{\D_{12}}{R \, T}c_1 \, \grad \mu_1.
\end{equation}
The chemical potential of $A_1$, say, is of the form $\mu_1 = \mu^0_1 + R\, T \ln (\gamma_1 x_1)$
with a reference chemical potential $\mu_1^0$ which only depends on pressure and temperature and
the so-called activity coefficient $\gamma_1 = \gamma_1 (x_1)$; note that the additional variable $c_{\rm tot}$
of  $\gamma_1$ is constant in the considered isobaric case. This yields
\begin{equation}
\label{binary3}
J_1= -\D_{12}  \Big( 1+\frac{x_1 \, \gamma_1' (x_1)}{\gamma_1 (x_1)} \Big) \, \nabla  c_1,
\end{equation}
where $\D_{12}$ is a function of $x_1$.
Inserting this into \eqref{mass-balance} leads to the nonlinear diffusion equation
\begin{equation}
\label{filtration}
\partial_t c_1 - \Delta \phi (c_1) =f(c_1),
\end{equation}
where the function $\phi :\R \to \R$ satisfies $\phi' (s c_{\rm tot})=\D_{12}(s) (1+s \gamma' (s)/\gamma (s))$ and, say, $\phi (0)=0$.
Equation (\ref{filtration}) is also known as the filtration equation (or, the generalized
porous medium equation) in other applications.
Note that \eqref{filtration} is locally wellposed in $L^1 (\Omega)$ as soon as $\phi$ is continuous and nondecreasing
which will also by used below; cf., e.g., \cite{Vazquez}.
For constant $\D_{12}$, the monotonicity of $\phi$ holds if $s\to s\gamma(s)$ is increasing.
This means that the chemical potential $\mu_1$ should be a monotone (increasing) function of $x_1$,
which characterizes systems without spontaneous phase separation.

A \emph{dilute mixture} is a mixture in which one component, say $A_P$, satisfies $x_P \approx 1$ and acts as a solvent, while
the other components are solutes and only appear in small concentrations, i.e.\ $x_i \ll 1$ for $i=1,\ldots , P-1$.
In this case the chemical potential of the dilute species is given by
\[
\mu_i = \mu_i^0 + RT \ln x_i.
\]
This leads to the diffusive fluxes
\begin{equation}
\label{dilute-flux}
{ J}_i=-\frac{\D_{iP}}{R \, T} \, c_{\rm tot}\, \nabla x_i =-\frac{\D_{iP}}{R \, T} \,\nabla c_i.
\end{equation}
Here the basic assumption is that interactions only occur between individual solutes and the solvent, but not between different solutes.
Hence, concerning the mixture composition, $\D_{iP}$ depends only on $x_i$ and $x_P$. Since $x_P$ is almost constant equal 1, it is essentially
a function of $x_i$, i.e.\ of $c_i$. This leads to Fick's law with diffusivities $d_i = d_i (c_i)$.

Combining the above prototype cases leads to a large class of mixtures in which two components are present in large amounts, while all other components are dilute.
This applies to many concrete cases in Chemical Reaction Engineering, in which
one species (e.g., water) acts as a solvent, one further species is the main feed into the process and the other constituents are further reactants, catalysts, initiators, intermediates or products.
This case leads to diffusivities $d_i =d_i (c_i, c_j)$ which not
only depend on $c_i$, but also on at least one further $c_j$, while still no cross-diffusion appears.

Let us note that other chemical applications as well as completely different motivations
also lead to RD-systems with concentration-dependent diffusivities.
Besides reactive turbulent flows (cf.\ \cite{BaldygaBourne}), let us only mention reactive transport
in the underground, i.e.\ inside porous media (cf.\ \cite{transport-in-pm}).
A common approach to model multicomponent transport in porous media employs an extension of the Maxwell-Stefan equations,
the so-called dusty gas model. The latter is based on adding another species, modeling the pore walls, which is immobile.
For a dilute species in a porous medium this again leads to diffusivities of type $d_i (c_i)$, as sketched above.

More general, system \eqref{mass-balance} can represent a set of population balances (cf., e.g., \cite{Murray}),
in which case $c_i$ denotes a number density of individuals of the $i$-th population.
Then the diffusive fluxes correspond to stochastic motions of the individuals, while additional
migrative fluxes might also occur in such situations.
Again, the $d_i$ will be non-constant as well as nonnegative.

Finally, it dependencies of the diffusion coefficients on other variables (like on the temperature)
are of relevance, but the evolution of these variables is not explicitly described (like in models
without energy balance to account for temperature changes), such dependencies can, in principle,
be incorporated via diffusivities $d_i = d_i (t,x,c)$ which depend also on $(t,x)$. This is
in particular used to include seasonal effects.\\[2ex]
\indent
The mass production terms $f_i$ involve the rate functions $r_i$ which
are nonlinear functions of the composition with superlinear growth, except in rare cases like for isomerizations of type $A_1\rightleftharpoons A_2$. Hence, while local-in-time existence of even classical solutions usually follows
from known results on quasi-linear parabolic PDE-systems (like the theory from \cite{amann89}, \cite{amann93}), the issue of global existence
of solutions can be a much more difficult one, depending on the structure of the reaction terms.
To this end, in order to have reliable information about the form of the $r_i$ at all, we focus on the case of (networks of)
elementary reactions. These are chemical reactions which run in a single step without the formation of intermediate
species. In other words, if intermediate steps occur, they have to be fully modeled by an appropriate reaction network.
In this case the rate functions for the elementary reactions are accurately modeled by so-called mass action kinetics.
To be more specific, the rate function $r$ for the single reversible reaction of type
\[
\alpha_1 A_1 + \ldots + \alpha_P A_P \rightleftharpoons \beta_1 A_1 + \ldots + \beta_P A_P
\]
with stoichiometric coefficients $\alpha_i, \beta_i \in \N_0$
is given as $r=r^f  - r^b$ with the forward and backward rates
\[
r^f (c)= k^f \prod_{i=1}^P c_i^{\alpha_i}  \qquad \mbox{ and } \qquad
r^b (c)= k^b \prod_{i=1}^P c_i^{\beta_i},
\]
respectively. It is important to note that throughout this paper, the so-called irreversible limit
case is also included, i.e.\ we allow for $k^f=0$ or $k^b=0$.

In case of $R$ chemical reactions, the stoichiometric coefficients are denoted
$\alpha_i^j, \beta_i^j$, where $i=1,\ldots ,P$ and $j=1,\ldots ,R$.
Correspondingly, the rate functions above become $r^f_j, r^b_j$. Then, the balance of molar
mass for species $A_i$ reads as
\[
\partial_t c_i + \div J_i = \sum_{j=1}^R (\beta_i^j - \alpha_i^j)(r^f_j - r^b_j).
\]
Notice that usually many of the coefficients $\alpha_i^j, \beta_i^j $ are in fact zero;
this is why they belong to $\N_0 = \N \cup \{0\}$.
The so-called rate constants $k^f_j, \, k^b_j$ are not constant but depend especially on the temperature.
Still, considering only isothermal systems, we will assume them to be constants below.

RD-systems with mass action kinetics, or more general rate functions of polynomial type, say,
but with {\it constant} Fickian diffusivities have been studied in many papers for long time.
Concerning global existence of solutions, already for constant diffusion coefficients the situation is complicated
unless all $d_i$'s are the same. A recent survey about the subject can be found in \cite{pierre10}.
Here, let us only emphasize that the main elementary reactions which occur in chemical reaction networks
are of the form
\begin{equation}
\label{A+B=P}
A_1 + A_2 \rightleftharpoons A_3
\end{equation}
or
\begin{equation}
\label{A+B=P+Q}
A_1 + A_2 \rightleftharpoons A_3 + A_4,
\end{equation}
i.e.\ at most two reaction partners appear on each side since (reactive) collisions of more than two molecules are very rare events (at least for moderate concentrations).
Note that we leave out reactions of the form $A_1 \rightleftharpoons A_2$ which are considered trivial due to
their linear rate functions, while $A_1 = A_2$ or $A_3=A_4$ is allowed in the reaction mechanism \eqref{A+B=P}, respectively
\eqref{A+B=P+Q}.
Reactions of type \eqref{A+B=P} occur for example if double bonds are opened in halogenizations, hydrations, sulfonizations etc., while
mechanism \eqref{A+B=P+Q} is typical for exchange reactions, where one reactant breaks into two parts, one of them being replaced
by the reaction partner.
Let us also note that a reaction which is formally of type \eqref{A+B=P+Q} might involve an intermediate species $A_5$, such that
the elementary steps are rather
\begin{equation}
\label{A+B-I-P+Q}
A_1 + A_2 \rightleftharpoons A_5 \rightleftharpoons A_3 + A_4,
\end{equation}
instead. In this case, the reaction is build from blocks of type \eqref{A+B=P}.
In fact, even without occurrence of an intermediate form $A_5$, 
the reaction from $A_1 + A_2$ to $A_3 + A_4$
proceeds via a so-called transition state, but the latter has a very limited life time of about $10^{-13} s$, only.
Compared to any transport process by diffusion, the transition hence is so fast that the transition state need not be
separately accounted for in the model. Indeed, the rigorous limit of the RD-system modeling \eqref{A+B-I-P+Q} as the intermediate's
life time approaches zero turns out to be the RD-system for \eqref{A+B=P+Q}; cf.\ \cite{BP-intermediate}.

For more information about chemical kinetics and reaction mechanisms see \cite{Espenson}.

As indicated above, the type of models considered here require small or moderate concentrations,
i.e.\ they loose their validity for large concentrations, especially in case of a blow-up.
Therefore, the question whether these models have the right intrinsic structure to prevent
solutions from blowing up, is a very natural one. This is of course equivalent to the
question of global existence of solutions, which in turn comes down to find
$L^\infty$-bounds for local solutions on any bounded time interval.

Global existence of solutions is known for a single reaction of type \eqref{A+B=P} in the case of {\it constant} diffusivities. Indeed,
it was shown in \cite{rothe} that for bounded initial data and space dimensions $N\leq 5$, the system \eqref{mass-balance}, \eqref{Neumann-BC}, \eqref{ICs} has a unique nonnegative classical solution, which is uniformly bounded. Global existence and boundedness in any space dimension for smooth $\Omega$ (of class $C^{2+\alpha}$, $0<\alpha<1$) and smooth initial data has been shown in \cite{feng91}. Both these approaches are based on semigroup theory and hence exploit the semilinear structure. This prototype RD-system also
has a particular triangular structure for which global existence of strong solutions is proved in \cite{pierre10} for more general systems, for any space dimension and bounded initial data.
This approach uses maximal $L^p$-regularity theory (see \cite{denk}) on the dual equations, and strongly relies on the linearity of the diffusion operators.

For a single reaction of type \eqref{A+B=P+Q}, still with constant diffusion coefficients,
the question of global existence of solutions has an affirmative answer only for $N=2$ so far,
while the physically more interesting case $N\geq 3$ is open; see 
\cite{GoudonVasseur}, where also the Hausdorff dimension of the set of possible singularities is estimated.

For {\it non-constant} diffusivities, the issue of global existence for such RD-systems is widely open.
The only closely related result which we are aware of is \cite{MorganWaggonner}, where the case $d_i(c_i)$ and reaction networks with
at most quadratic terms and an appropriate triangular-type structure ("intermediate sum"-condition) are considered and global existence is obtained in case $N=2$.

In the present paper, we consider reaction networks with building blocks of type \eqref{A+B=P} and with diffusivities which depend
on time, space and composition. We obtain global existence of solutions for initial values from an appropriate Sobolev space, the
regularity index of which is optimal in a certain sense. In case of bounded diffusion coefficients of type $d_i (c_i)$ and $N=3$,
we can relax to arbitrary $L^\infty$-bounded initial values. The admissible space dimensions are always at least $N=3$.
The core point of our approach is a thorough analysis of the RD-system with a single reversible reaction of type \eqref{A+B=P}. We first derive an initial estimate on the solutions from the conservation of the total mass for general diffusivities, and from $L^2$-techniques in the case of diffusivities $d_i(c_i)$. Since the solutions are nonnegative and the reaction terms for some equations are linearly bounded from above, by classical results on parabolic equations, this initial estimate may be improved for the corresponding $c_i$. For small space dimensions, this provides new estimates on some quadratic reaction terms, which allows to improve the regularity for other concentrations until an increased regularity is obtained for all species.
Bootstrapping this procedure, we may estimate the solution in $L^p((0,T) \times \Omega  )$ for any $T>0$, $p<+\infty$, and then in $L^\infty((0,T)\times \Omega )$; cf. \cite{LSU}. Global existence then follows from a global existence criterion due to \cite{amann93}.
%

Finally, let us mention some related works, which all concern the case of constant diffusivities: asymptotics has first been studied by Rothe in \cite{rothe}, where it is proved that $c(t)$ converges to a uniquely determined homogeneous stationary state when $t\rightarrow +\infty$. In \cite{DF06}, Desvillettes and Fellner used the entropy method to give explicit convergence rates to the equilibrium. The present paper does not employ entropy-like Lyapunov functionals and, hence, also
applies to irreversible limits.

The fast-reaction limit $k^f,k^b\rightarrow +\infty $ for the RD-system \eqref{syst} has first been studied in \cite{BoHabil}, in the special case when the diffusion coefficients are equal, and then in \cite{BPR12} for the case of different but constant diffusivities. Note in passing that the techniques developed in the latter paper carry over with only slight modifications to the case of nonlinear diffusions of the type $d_i(c_i)\nabla c_i$. Then using the above global existence result, Theorem 1 in \cite{BPR12} can be extended to the case of diffusivities $(\ref{hypd_i2})$ for space dimensions $N\leq9$.

In Section \ref{S2}, we prove well-posedness for the reaction-diffusion system associated with the reaction
$A_1 + A_2 \rightleftharpoons A_3$ for certain space dimensions including $N\geq 3$. For clarity of presentation, some technical details are postponed to an appendix.

Section \ref{S3} is devoted to the case of reaction networks involving $P$ chemically reacting species $A_1,\ldots,A_P$, where chemical reactions are assumed to be of the type $A_i+A_j\rightleftharpoons A_k$ and the total mass of involved atoms is preserved. After re-sorting the reactions and chemical species to get a block-triangular structure, we prove that the ideas developed in Section \ref{S2} can be adapted to this case, but under stronger restrictions on the admissible space dimensions.

In Section~\ref{sec:ext} we focus on the quasi-linear case $d_i(c_i)$ and show that, for $N=3$, Theorem~\ref{theomain} and Theorem~\ref{th2}
stay valid for initial values being merely bounded and measurable. Moreover, again focusing on $N=3$,
we extend to rate functions of type $f(c_1)g(c_2)-h(c_3)$, satisfying growth conditions of fractional
power type with appropriate conditions for the exponents.

Finally, the Appendix contains the proof of Lemma \ref{lemEstimates} and the proof of a theorem from \cite{LSU}, but adapted to Neumann boundary conditions.

\section{The Rothe system with variable diffusivities}\label{S2}
\subsection{Global well-posedness for quasi-linear and semi-linear systems}
\nd We consider the reaction-diffusion system
\begin{eqnarray}\label{syst}
\left\{
\begin{array}{l}
\left.
\begin{array}{l}
\partial_t c_1-\mathrm{div}(d_1(t,x,c)\nabla c_1)=-k^f c_1c_2+ k^b c_3\\[0.5ex]
\partial_t c_2-\mathrm{div}(d_2(t,x,c)\nabla c_2)=-k^f c_1c_2+ k^b c_3\\[0.5ex]
\partial_t c_3-\mathrm{div}(d_3(t,x,c)\nabla c_3)=+k^f c_1c_2- k^b c_3\\[0.5ex]
\end{array}
\right\} \vspace{0.1cm}\; \mbox{ on }\; (0,+\infty)\times\Omega,\\
\displaystyle
\partial_\nu c_1=\partial_\nu c_2=\partial_\nu c_3=0\;\mbox{ on }\;(0,+\infty)\times\partial\Omega,\\[1.5ex]
c=(c_1,c_2,c_3);\ c(0,\cdot)=(c_{0,1},c_{0,2},c_{0,3}) \on \Omega,\ c_{0,i}\geq 0.
\end{array}
\right.
\end{eqnarray}
Throughout the paper, $\Omega$ denotes an open and bounded subset of $\R^N$, whose boundary $\partial \Omega$ is supposed to be at least of class $C^2$. The normal exterior derivative of a function $c$ on $\partial\Omega$ is denoted by $\partial_\nu c$.
As mentioned in the introduction, the system $(\ref{syst})$ represents the time-evolution of the concentration $c=(c_1,c_2,c_3)$ of three chemical species taking part in the reaction
$$A_1+A_2 \underset{k^b}{\overset{k^f}{\rightleftharpoons}} A_3,$$
where $k^f,k^b\geq 0$ are the rate constants for the forward and backward reaction.
Recall that $k^f=0$ or $k^b =0$ refers to the irreversible limit case which is included.
The reaction rates are modeled by mass action kinetics, which is usually relevant for such an elementary reaction. The transport of species is assumed to be driven only by diffusion, with mass fluxes of the type $d_i(t,x,c)\nabla c_i$. Observe that indirect cross-effects can occur, since the diffusion coefficients depend on all species. This simple system is interesting since it contains most mathematical difficulties to treat the case of larger networks of reactions of the type $A_i+A_j\rightleftharpoons A_k$, satisfying atomic conservation (see Section~\ref{S3}).\vspace{2mm}

The aim of this work is to prove the well-posedness of system $(\ref{syst})$ for nonlinear diffusivities and smooth initial data. More precisely, we assume that the diffusion coefficient for the $i^{th}$ species $d_i=d_i(t,x,c)$ depends on all the concentrations and
\begin{equation}\label{hypd_i1}
d_i \in C^{2-} ([0,+\infty)\times \overline \Omega \times \R^3,\R_+)\;\mbox{; }\;\exists \underline d>0 \mbox{ such that } \underline{d}\leq d_i,
\end{equation}
where for $k\geq 1$, $ C^{k-}$ is the space of $(k-1)$ times continuously differentiable functions whose derivatives of order $k-1$ are Lipschitz continuous. The special situation when $d_i$ only depends on the $i^{th}$ variable (i.e.\ $d_i=d_i(c_i)$) is particularly interesting since it allows to use some recent $L^2$-techniques, which are not available in general. In this case, we  write $d_i(c_i)$ instead of $d_i(t,x,c)$ and assume
\begin{equation}\label{hypd_i2}
d_i \in C^{2-} ( \R,\R_+)\;\mbox{; }\;\exists \underline d>0 \mbox{ such that } \underline{d}\leq d_i.
\end{equation}
The first step in the proof is of course the local existence of solutions which is based on a local well-posedness result from Amann \cite{amann93}, where the following notion of weak solution is used: consider the general reaction-diffusion system
\begin{eqnarray}\label{systintro}
\left\{
\begin{array}{rcll}
\partial_t c_i-\mathrm{div}(d_i(t,x,c)\nabla c_i)&=&f_i(c)\vspace{0.1cm}\;	 &\mbox{ on }\; (0,+\infty)\times\Omega,\\
\partial_\nu c_i&=& 0\;								&\mbox{ on }\;(0,+\infty)\times\partial\Omega,\\[0.5ex]
c_i(0,\cdot)&=&c_{0,i}								&\on \Omega,
\end{array}
\right.
\end{eqnarray}
where $c=(c_1,\ldots,c_P)$ and $f_i\in C^{1-}(\R^P)$.\saut
{\nd \bf Definition (weak $ \mathbf{W^s_p}$-solution). }{\it
Let $T\in (0,+\infty],\;p> 1,\; p'=p/(p-1),\;s>0$ satisfying
\begin{equation}\label{intro1}
 \frac{N}{p}<s<\min(1+\frac{1}{p},2-\frac{N}{p})\;,
\end{equation}
and assume $c_{0,i} \in W_p^s(\Omega)$. A weak $W^s_p$-solution of system $(\ref{systintro})$ on $[0,T)$ is a function $c=(c_1,\ldots,c_P):[0,T)\times \Omega \rightarrow \R^P$ such that
\begin{equation*}
 c\in {C}([0,T);W^s_p(\Omega)^P)\cap {C}^1((0,T); W^{s-2}_p(\Omega)^P),
\end{equation*}
$c(0)=c_0$ and for all $t\in (0,T)$, $v\in W^{2-s}_{p'}(\Omega)$, $i\in \{1,\ldots,P \}$,
\begin{equation*}
 \langle \partial_t c_i(t),v\rangle_{W_p^{s-2},W_{p'}^{2-s}} +\langle d_i(t,x,c)\nabla c_i(t),\nabla v\rangle_{W_p^{s-1},W_{p'}^{1-s}} =\langle f_i(c),v\rangle_{L^\infty,W_{p'}^{2-s}}.
\end{equation*}
}\\[-2ex]
The above definition requires slightly more regularity than the corresponding notion in \cite{amann93} in order to avoid using a third function space in defining the regularity class.
But solutions of the quasilinear PDEs considered in \cite{amann93} are in fact more regular and,
in particular, fulfil the assumptions of our definition.

Throughout the rest of the paper, by a {\bf classical solution} we denote a function that belongs to
$C ([0,T)\times \overline \Omega)\cap C^1((0,T); C(\overline \Omega))\cap C((0,T); C^2(\overline \Omega))$ and satisfies the equations pointwise.\\

Since we are interested in systems describing chemical concentrations, the nonnegativity of the solutions has to be preserved. For the general system $(\ref{systintro})$, it is well-kown that a necessary and sufficient condition is the so-called quasi-positivity of the reaction terms: \saut
{\nd \bf Definition (quasi-positivity). }{\it
 A vector field $f=(f_1,\ldots,f_P):\R^P\rightarrow \R^P,\; y=(y_1,\ldots,y_P)\mapsto f(y)$ is {\it quasi-positive} if
 \begin{equation}\label{quasipositivity}
  \forall i\in \{1,\ldots,P\},\forall y\in \R_+^P,\quad y_i=0 \Rightarrow f_i(y)\geq 0.
 \end{equation}
}
$\mbox{ }$\\[-1ex]
\noindent
We can now state the first main result.
Throughout this paper, for any $T>0$, we use the common notation $Q_T=(0,T)\times \Omega$ and $\Sigma_T=(0,T)\times \partial \Omega$.
\begin{theo}\label{theomain} Let $p> 1$, $s>0$ satisfying $(\ref{intro1})$ and $c_0 \in W_p^s(\Omega,\R_+^3)$. System $(\ref{syst})$ has a unique global weak-$W^s_p$ solution $c=(c_1,c_2,c_3):[0,+\infty)\times \Omega \rightarrow \R^3$ provided one of the following conditions is satisfied:
\begin{enumerate}[\quad (i)]
\item  $N\leq 5$ and the diffusivities $d_i(t,x,c)$ satisfy $(\ref{hypd_i1})$.
\item $N\leq 9$ and the diffusivities $d_i(c_i)$ satisfy $(\ref{hypd_i2})$.
\end{enumerate}
This solution is nonnegative. It is actually a classical solution and $(\ref{syst})$ is satisfied in a pointwise sense. Moreover
\begin{equation}\label{Linf:rothe}
 \forall T>0,\; \exists C=C(\|c_0\|_{L^\infty(\Omega)^3},T )>0 \;\mbox{ such that }\;\|c\|_{L^\infty (Q_T)^3}\leq C.
\end{equation}
If, in addition, $d_i$ and $\partial \Omega$ are smooth, then  $c\in {C}^\infty((0,+\infty)\times \overline\Omega\;;\R_+^3)$.
\end{theo}
{\nd\it Outline of the proof.} According to Amann's theory \cite{amann93}, local well-posedness and nonnegativity holds for $(\ref{syst})$. The solution is global provided it is  bounded in $L^\infty(Q_T)$ for any $T<+\infty$. The conservation of the total mass gives a first estimate on $c$ in $L^\infty(0,T;L^1(\Omega)^P)$, and actually the reaction terms in $(\ref{syst})$ are bounded in $L^1(Q_T)$.
Then we use the theory of scalar parabolic equations to estimate $c$ in $L^{{(N+2)}/{N}-\eps}(Q_T)$ for any $\eps>0$. Reaction terms for $c_1$ and $c_2$ are (linearly) bounded above by $c_3$, so $c_1$ and $c_2$ can be estimated in a better $L^p(Q_T)$-space ($p$ depending on $N$). Then the reaction term for $c_3$ is bounded above by $c_1c_2$, and for small enough space dimensions the previous estimates are sufficient to improve the regularity on $c_3$. Bootstrapping this procedure, we get estimates in $L^p(Q_T)$ for any $p$ and whence in $L^\infty(Q_T)$ for any $T>0$ by classical results
from \cite{LSU}. In the special case of diffusivities $d_i(c_i)$, we can directly start with estimates in $L^2(Q_T)$ which allows for higher space dimensions.

\subsection{Proof of Theorem \ref{theomain}}
\nd Notice first that, using the rescaling
$$(t,x)\mapsto \frac{k^f}{k^b}c(\frac{t}{k^b},x),$$
we can assume, without loss of generality, that $k^f=k^b=1$. As mentioned above, the reaction term in $(\ref{syst})$ satisfies the quasi-positivity assumption $(\ref{quasipositivity})$, so according to Amann's theory (see \cite{amann93}, Theorems 14.4 and 15.1, \cite{amann89} for the proofs), $(\ref{syst})$ has a unique nonnegative weak $W^s_p$-solution $c$, defined on a maximum time interval $[0,T^*)$, $T^*\leq +\infty$. The additional regularity properties as
stated in Theorem \ref{theomain} are consequences of Theorem 14.6 and Corollary 14.7 in \cite{amann93}.\saut
It remains to prove that the solution is global and, according to Theorem 16.3 in \cite{amann93}, it suffices to prove that $c$ is  bounded in $L^\infty(Q_T)^3$ for any $T>0$. For this purpose, we first estimate the solution in $L^p(Q_T)$ spaces for finite $p$. The subsequent Lemma is the main tool to improve these estimates by a bootstrap procedure: given a bound in $L^r(Q_T)$ on the positive part of a reaction term $f_i$, it shows in which $L^q(Q_T)$ space $c_i$ is bounded.
In the case $r>1$, the resulting regularity corresponds to the maximal $L^p$-regularity, possibly
up to inclusion of the limit case.
The proof is given in the Appendix.
\begin{lem}\label{lemEstimates}
Let $d:Q_T\rightarrow \R_+$ be measurable and such that $\underline d\leq d$ for some $\underline d>0$. Let $f\in L^{r}(Q_T)$ for $1\leq r<+\infty$ and $u$ be a nonnegative classical solution of
\begin{equation}\label{eqlemEstimates}
\partial_t u -\mathrm{{div}}(d(t,x)\nabla u)\leq f (t,x) \mbox{ in } Q_T\;;\quad
\partial_\nu u =0 \on \Sigma_T\;;\quad
 u(0)=u_0\in L^\infty(\Omega).
\end{equation}
Then $\|u\|_{L^{q}(Q_{T})}$ is bounded by a constant depending only on $T,\underline{d},\|f\|_{L^{r}(Q_{T})}$ and $ \|u_0\|_{L^\infty(\Omega)}$, provided $1\leq q<+\infty$ and $(r,q)$ satisfies one of the following conditions:
$$
(i)\;\;r=1 \mbox{ and }
\left\{
\begin{array}{ll}
 1-\frac{2}{N+2}<\frac{1}{q} 	& \mbox{ for }N \geq 2\\[0.8ex]
q<2				&\mbox{ for }N=1
\end{array}
\right.
\;;\quad (ii)\;\; r>1 \mbox{ and }
\left\{
\begin{array}{ll}
 \frac{1}{r}-\frac{2}{N+2}\leq \frac{1}{q} 			& \mbox{ for }N \geq 3\\[0.8ex]
\frac{1}{r}-\frac{1}{2}<\frac{1}{q}				& \mbox{ for }N=2\\[0.8ex]
\frac{1}{r}-\frac{1}{2}\leq \frac{1}{q}				&\mbox{ for }N=1
\end{array}
\right. .
$$
\end{lem}
{\nd \bf Step 1. The initial estimate.}\vspace{1mm}\\
Given $0<T<+\infty$ with $T\leq T^*$, we estimate $c$ on $Q_T$ as follows.\vspace{1mm}\\
{\it For diffusivities $d_i(t,x,c)$ satisfying} $(\ref{hypd_i1})$:\vspace{1mm}\\
Let $r_0\in [1,(N+2)/N)$ if $N\geq 2$, $r_0\in [1,2)$ if $N=1$, and let us prove that
\begin{equation} \label{step1}
 \exists C=C(T,\underline d,\|c_0\|_{L^\infty(\Omega)^3})>0\;:\quad \|c\|_{L^{r_0}(Q_T)^3}\leq C.
\end{equation}
Using the homogeneous Neumann boundary conditions in $(\ref{syst})$, it is clear that
$$\frac{d}{dt}\int_\Omega \big( c_1(t)+c_2(t)+2c_3(t)\big) =0.$$
As $c$ is nonnegative,
\begin{equation}\label{step1L1bound}
\forall i\in \{1,2,3\},\quad \sup_{t\in [0,T^*)}\|c_i(t)\|_{L^1(\Omega)}\leq \|c_{0,1}\|_{L^1(\Omega)}+\|c_{0,2}\|_{L^1(\Omega)}+2\|c_{0,3}\|_{L^1(\Omega)}.
\end{equation}
After integration of the first equation in $(\ref{syst})$ on $Q_T$ and integration by parts,
\begin{equation*}
 \int_{Q_T}c_1c_2=\int_{Q_T}c_3 +\int_\Omega c_{0,1}-\int_\Omega c_1(T).
\end{equation*}
All the integrals on the right-hand side are bounded, so $c_1c_2$ is bounded in $L^1(Q_T)$, and the reaction terms in $(\ref{syst})$ are bounded in $L^1(Q_T)$. Then $(\ref{step1})$ is a consequence of Lemma \ref{lemEstimates} $(i)$.\saut
{\nd \it With diffusivities $d_i(c_i)$ satisfying} $(\ref{hypd_i2})$:\saut
Let us prove that
\begin{equation}\label{step1bis}
 \exists C=C(T,\underline d,\|c_0\|_{L^2(\Omega)^3})>0:\quad \|c\|_{L^2(Q_T)^3}\leq C.
\end{equation}
In this case, $(\ref{syst})$ can be rewritten as
\begin{equation}\label{systdici}
\partial_t c_i -\Delta D_i(c_i)=\zeta_i(c_1c_2-c_3)\on Q_T\;;\quad \partial_\nu D_i(c_i) =0 \on \Sigma_T\;;\quad c_i(0)=c_{0,i}\on \Omega,
\end{equation}
where $i\in \{1,2,3\}$, $\zeta=(-1,-1,1)$, $D_i(y)=\int_0^y d_i(s)ds$. Using assumption $(\ref{hypd_i2})$, $\underline{d} \, y\leq D_i(y)$ for $y\geq0$. Then $(\ref{step1bis})$ is a straightforward consequence of the following lemma (applied to $(c_1,c_2,2c_3)$), which generalizes Proposition 6.1 in \cite{pierre10} to the case of nonlinear diffusion of filtration equation type:
\begin{lem}\label{lemL2}
Let $T>0$, $c=(c_1,\ldots,c_P)$ be a nonnegative solution of
\begin{equation}\label{eqlemL2}
 \partial_tc_i-\Delta D_i(c_i)=f_i\on Q_T\;;\quad \partial_\nu D_i(c_i)=0 \on \Sigma_T\;;\quad c_i(0)=c_{0,i}\in L^2(\Omega,\R_+),
\end{equation}
where $i\in \{1,\ldots,P\}$, $f_i:Q_T\rightarrow \R$ is measurable
with $\sum_{i=1}^P f_i\in L^2(Q_T)$, $D_i:\R_+\rightarrow \R_+$ and
\begin{equation}
\exists \underline d >0:\;\forall y\geq 0,\;\underline d \,y\leq D_i(y).
\end{equation}
Then there exists $C=C(T,\underline d,\|\sum_{i=1}^P f_i\|_{L^2(Q_T)},\|c_0\|_{L^2(\Omega)^P})>0$ such that
\begin{equation}\label{eqL2estimate}
 \|c\|_{L^2(Q_T)^P}\leq C.
\end{equation}
\end{lem}
{\nd\bf Proof of Lemma \ref{lemL2}. }Set
$$  W:=\sum_{i=1}^P c_i \;;\quad W^0:=\sum_{i=1}^ P c_{i}^0\;;\quad A:=\frac{\sum_{i=1}^P D_i(c_i)}{\sum_{i=1}^P c_i}\;;\quad F:=\sum_{i=1}^P f_i$$
and note that $A\geq \underline d$.
Let $t\in (0,T)$ and integrate $(\ref{eqlemL2})$ on $(0,t)$ to get, for $i\in \{1,\ldots,P\}$,
\begin{equation}
c_i-\Delta \int_0^t D_i(c_i)=c_{i}^0+\int_0^tf_i \on Q_T\;;\quad \partial_\nu D_i(c_i)=0 \on \Sigma_T\;;\quad c_i(0)=c_{0,i} \on \Omega.
\end{equation}
Summing these equations over $i$ yields
\begin{equation}
 W-\Delta \int_0^t AW =W^0+\int_0^t F\on Q_T\;;\quad \partial_\nu (AW)=0\on \Sigma_T\;;\quad W(0)=W^0\on \Omega.
\end{equation}
After multiplication by $AW$, integration over $Q_T$ and integration by parts, we get
\begin{align}\nonumber
 \int_{Q_T} AW^2+\int_{Q_T}\nabla(AW) \cdot \nabla \int_0^t AW &=\int_{Q_T} W^0AW+\int_{Q_T}\left(\int_0^t F\right)AW,\nonumber
\end{align}
which yields
\begin{align}\nonumber
\int_{Q_T} AW^2+\frac{1}{2}\int_{\Omega}\left|\nabla \int_0^T AW\right|^2 &=\int_{\Omega} W^0\int_0^TAW+\int_{Q_T} F\int_t^T AW\nonumber\\
					&\leq \|W^0\|_{L^2(\Omega)} \|\int_0^T AW\|_{L^2(\Omega)} +\sqrt{T}\|F\|_{L^2(Q_T)}\|\int_0^TAW\|_{L^2(\Omega)}\nonumber\\
					&\leq C \|\int_0^T AW\|_{L^2(\Omega)},
\end{align}
where $C>0$ denotes a constant depending only on $\|F\|_{L^2(Q_T)},\|c_0\|_{L^2(\Omega)^P}$, $\underline d $ and $T$. Using the Poincar\'e-Wirtinger inequality,
\begin{align*}
 \exists C>0:\quad \int_{Q_T} AW^2+\frac{1}{2}\int_{\Omega}\left|\nabla \int_0^T AW\right|^2\leq C\left(\|\nabla \int_0^T AW\|_{L^2(\Omega)}+\int_{Q_T}AW \right)\nonumber.
\end{align*}
Then Young's inequality yields
\begin{align}
 \exists C>0:\quad \int_{Q_T} AW^2+\frac{1}{2}\int_{\Omega}\left|\nabla \int_0^T AW\right|^2\leq C+\frac{1}{4}\int_{\Omega}\left|\nabla \int_0^T AW\right|^2+C\int_{Q_T}AW \label{eqLemL2a}.
\end{align}
Letting $\alpha>0$, $ \{W> \alpha\}:=\{(t,x)\in Q_T:\; W(t,x)> \alpha\}$ and $\{W\leq  \alpha\}:=Q_T\backslash \{W> \alpha\}$, we have
\begin{align}
 \int_{Q_T}AW&=\int_{ \{W> \alpha\}}AW +\int_{\{W\leq  \alpha\}}AW\nonumber\\
		&\leq \frac{1}{\alpha}\int_{Q_T}AW^2 +\int_{\{W\leq  \alpha\}}\sum_{i=1}^pD_i(c_i)\nonumber\\
		&\leq \frac{1}{\alpha}\int_{Q_T}AW^2+M_\alpha\label{eqLemL2b},
\end{align}
where we used $c_i\leq \alpha$ on $\{W\leq  \alpha\}$ and set
$M_\alpha :=|\Omega| \, T \max_{0\leq s\leq \alpha}\sum_{i=1}^P D_i(s)$.
Choosing $\alpha =2C	$, where $C$ is the constant from $(\ref{eqLemL2a})$, we get
$$ \underline d\int_{Q_T}W^2\leq \int_{Q_T} AW^2+\frac{1}{2}\int_{\Omega}\left|\nabla \int_0^T AW\right|^2\leq 2C(M_{2C}+1). $$
Using $c_i\geq 0$ and $W=\sum_{i=1}^P c_i$, this proves the desired bound on $c$ in $L^2(Q_T)^P$.\cqfd

{\nd \bf Step 2. The bootstrap procedure.}\vspace{1mm}\\
Let us prove that the maximal solution of $(\ref{syst})$ is bounded in $L^p(Q_T)$ for any $p<+\infty$ and any $T\leq T^*$, $T<+\infty$. The idea is to exploit the fact that the reaction terms for $c_1 $ and $c_2$ are linearly bounded from above to get new estimates on $c_1$ and $c_2$. For small space dimensions, we get a sufficiently strong estimate on $c_1c_2$, which is an upper bound for the reaction term for $c_3$, such that we can improve the estimate on $c_3$. Then we go back to the equations in $c_1$ and $c_2$ and bootstrap this procedure.\saut
Assume first that $N=1$. For diffusivities satisfying $(\ref{hypd_i1})$ or $(\ref{hypd_i2})$, according to $(\ref{step1})$, $c$ is bounded in $L^{r_0}(Q_T)^3$ for $r_0<2$. Using Lemma $\ref{lemEstimates}$, $c_1$ and $c_2$ are bounded in $L^{p}(Q_T)$ for any $p<+\infty$, so $c_1c_2$ is also bounded in any $L^{p}(Q_T)$ and, using once more Lemma $\ref{lemEstimates}$, $c_3$ is bounded in any $L^p(Q_T)$.\saut
For $N\geq 2$, let $r_0>1$ be such that $c$ is bounded in $L^{r_0}(Q_T)^3$. According to Lemma \ref{lemEstimates},
\begin{enumerate}[$\quad$]
\item $c_1,c_2$ are bounded in $L^{q_1}(Q_T)$, where $ \frac{1}{r_0}-\frac{2}{N+2}<\frac{1}{q_1}$;
\item $c_1c_2$ is bounded in $L^{q_2}(Q_T)$, where $ \frac{2}{r_0}-\frac{4}{N+2}<\frac{1}{q_2}$. We can choose $q_2\geq 1$ provided
\begin{equation}\label{cond1}
 \frac{2}{r_0}-\frac{4}{N+2}< 1;
\end{equation}
\item $c_3$ is bounded in $L^{r_1}(Q_T)$, where
\begin{equation}\label{eqimprovement}
\frac{2}{r_0}-\frac{6}{N+2}<\frac{1}{r_1}.
\end{equation}
\end{enumerate}
The initial estimate is improved if we can choose $r_0<r_1$, i.e.\ if
\begin{equation}\label{cond2}
\frac{1}{r_0}<\frac{6}{N+2}.
\end{equation}
Suppose $r_0$ satisfies conditions $(\ref{cond1})$ and $(\ref{cond2})$. Then $c$ is bounded in $L^{r_1}(Q_T)^3$ for some $r_1>r_0$, which also satisfies $(\ref{cond1})$ and $(\ref{cond2})$. Then it is clear that we can build by induction an increasing sequence $(r_n)_{n\in\N}$ such that $c$ is bounded in $L^{r_n}(Q_T)^3$ and
$$\frac{2}{r_n}-\frac{6}{N+2}<\frac{1}{r_{n+1}}.$$
Let us prove that $ (r_n)_{n\in\N}$ can be built such that $r_n\rightarrow +\infty$. Let $0<\eps <\frac{6}{N+2}-\frac{1}{r_0}$. We define $r_{n+1}>r_n$ by
\begin{enumerate}[\qquad]
\item If $\;\frac{2}{r_n}-\frac{6}{N+2}< 0$\;,\quad $r_{n+1}=r_n+1$.
\item If $\;\frac{2}{r_n}-\frac{6}{N+2}\geq0$\;, \quad$\frac{1}{r_{n+1}}=\frac{2}{r_n}-\frac{6}{N+2}+\eps$.
\end{enumerate}
Suppose that $\frac{2}{r_n}-\frac{6}{N+2}\geq 0$ for all $n\in\N$. Then $u_n:=\frac{1}{r_n}\in (0,1]$ is decreasing and satisfies $u_{n+1}=2u_n-\frac{6}{N+2}+\eps$. This yields $u_n \rightarrow -\infty$, a contradiction, so there exists $n_0\in\N$ such that $\frac{2}{r_{n_0}}-\frac{6}{N+2}<0$. Then for all $n\geq n_0$, $r_{n}=r_{n_0}+n-n_0$ and therefore $r_n\rightarrow +\infty$.
Consequently, $c$ is bounded in $L^p(Q_T)^3$ for any $p<+\infty$. \saut
It remains to give some explicit sufficient conditions so that we can choose $r_0$ satisfying $(\ref{cond1})$ and $(\ref{cond2})$:
\begin{enumerate}[(i)]
\item For diffusivities $d_i(t,x,c)$ satisfying $(\ref{hypd_i1})$: according to $(\ref{step1})$, $c$ is bounded in $L^{r_0}(Q_T)^P$ for $r_0<\frac{N+2}{N}$ (since $N\geq 2$). Hence equations $(\ref{cond1})$ and $(\ref{cond2})$ can be satisfied if and only if $N<6$.
\item For diffusivities $d_i(c_i)$ satisfying $(\ref{hypd_i2})$: according to $(\ref{step1bis})$, $c$ is bounded in  $L^{r_0}(Q_T)^P$ with $r_0=2$. Hence equations $(\ref{cond1})$ and $(\ref{cond2})$ are satisfied if and only if $N<10$.
\end{enumerate}
{\bf Step 3.} Once we know that $c$ is bounded in $L^p(Q_T)^P$ for any $p<+\infty$, we can use a classical result from \cite{LSU} on parabolic equations (see Theorem III.7.1 there
and Theorem~\ref{ch5:th} below) to say that for all $i$, $c_i$ is bounded in $L^\infty(Q_T)$. This is valid for any $T\leq T^*$, $T<+\infty$, so using Theorem 16.3 in \cite{amann93}, $T^*=+\infty$, i.e.\ $c$ is a global solution.\\[-2ex]
$\mbox{ } $\hfill $\Box$\\[1ex]

\begin{rem}
\begin{enumerate}
\item
In \cite{LSU}, Theorem III.7.1 is stated for Dirichlet boundary conditions. The result also holds for Neumann boundary conditions.
For completeness and self-containment of this paper, we include the proof for the latter case in the Appendix.
\item
In \cite{amann93}, the results we used from Chapters 14, 15 and 16 are stated for time-independent operators. To see that they are still valid for the time-dependent case, it is sufficient to ``artificially'' add the time in the equations, replacing $c=(c_1,\ldots,c_P)$ by $\tilde c=(c_1,\ldots,c_P,c_{P+1})$ in $(\ref{systintro})$, where $c_{P+1}$ satisfies
$\partial_t c_{P+1} -\Delta c_{P+1} =1$ with homogeneous Neumann boundary conditions
and $c_{P+1} (0)=0$. Note that $c_{P+1}(t,x)\equiv t$, then.
\item In the case of Michaelis-Menten-Henri (MMH) enzymatic reaction
$$A_1+A_2 \underset{k_{-1}}{\overset{k_1}{\rightleftharpoons}}A_3 \underset{k_{-2}}{\overset{k_2}{\rightleftharpoons}}A_1+A_4\;;\hspace{3mm} \tiny{ k_1,k_{-1},k_2,k_{-2}\geq 0},$$
we are led to the equations
\begin{eqnarray*}
\left\{
\begin{array}{l}
\left.
\begin{array}{lrccc}
\partial_t c_1-\mathrm{div}(d_1(t,x,c)\nabla c_1)=&-k_1c_1c_2	&+k_{-1}c_3	&+k_2c_3	 &-k_{-2}c_1c_4\\[0.5ex]
\partial_t c_2-\mathrm{div}(d_2(t,x,c)\nabla c_2)=&-k_1c_1c_2	&+k_{-1}c_3	&		 &\\[0.5ex]
\partial_t c_3-\mathrm{div}(d_3(t,x,c)\nabla c_3)=&k_1c_1c_2	&-k_{-1}c_3	&-k_2c_3	 &+k_{-2}c_1c_4\\[0.5ex]
\partial_t c_4-\mathrm{div}(d_4(t,x,c)\nabla c_4)=&	& 		 &+k_2c_3	 &-k_{-2}c_1c_4\\[0.5ex]
\end{array}
\right\} \vspace{0.1cm}\; \mbox{ on }\; (0,+\infty)\times\Omega,\\
\displaystyle
\partial_\nu c_1=\partial_\nu c_2=\partial_\nu c_3=\partial_\nu c_4=0\;\mbox{ on }\;(0,+\infty)\times\partial\Omega,\\[1.5ex]
c_i(0,\cdot)=c_{0,i},\ c_{0,i}\in L^\infty(\Omega,\R_+).
\end{array}
\right.
\end{eqnarray*}
Similarly as in $(\ref{syst})$, the reaction terms for $c_1$, $c_2$ and $c_4$ are linearly bounded above, and it is clear that with obvious modifications in the above proof, the results from Theorem \ref{theomain} also hold for this system, whith the same space dimension restrictions. In the literature on MMH reaction systems, the second reaction is usually assumed to be irreversible with $k_{-2}=0$. Note that this case is included in our analysis.
\end{enumerate}
\end{rem}

\section{Networks of elementary reactions}\label{S3}
In this section, we suppose that {$P$} chemical species $A_1,\ldots,A_P$ are present, and that they are involved in $R$ chemical reactions of the type
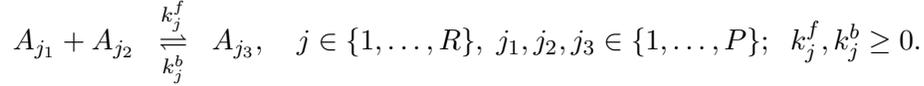
\begin{equation*}
A_{j_1}+A_{j_2} \;\;\underset{k^b_j}{\overset{k^f_j}{\rightleftharpoons}}\;\; A_{j_3},\quad j\in \{1,\ldots,R\},\; j_1,j_2,j_3\in\{1,\ldots,P\};\;\; k^f_j,k^b_j\geq 0.
\end{equation*}
Remark that $j_1$ and $j_2$ are not necessarily distinct, so that reactions of the type $2A_{j_1}\rightleftharpoons A_{j_3}$ are included, as well as the irreversible reactions $A_{j_1}+A_{j_2}\rightarrow A_{j_3}$ and $A_{j_3}\rightarrow A_{j_1}+A_{j_2}$, which are obtained by taking $k^b_j=0$, respectively $k^f_j=0$.

As before, $c_i$ denotes the concentration of species $A_i$. Let $(\eps_1,\ldots,\eps_P)$ be the canonical basis of $\R^P$ and define the so-called {\it stoichiometric vectors} as $\alpha_j:=\eps_{j_1}+\eps_{j_2}$, $\beta_j:= \eps_{j_3}$ and $\nu_j:=\beta_j-\alpha_j$. The {\it stoichiometric matrix }$M\in \R^{P\times R}$ is the matrix whose columns are $\nu_1,\ldots,\nu_R$. On the basis of mass action kinetics, the reaction rate for the $j^{th}$ reaction is given by $r_j(c)=k^f_j c_{j_1}c_{j_2}-k^b_j c_{j_3}$.\\
We also assume that an {\it atomic conservation law} (see \cite{erdi}, Chapter 3) applies:
we impose the condition
\begin{equation}\label{conservationlaw}
 \exists e\in (0,+\infty)^P:\; \forall i\in \{1,\ldots,R\}, \quad \langle e,\nu_i\rangle =0,
\end{equation}
where $\langle \cdot,\cdot\rangle$ denotes the scalar product on $\R^P$.
Note that assumption $(\ref{conservationlaw})$ excludes chemical reactions of the type $A_{j_1}+A_{j_2}\rightleftharpoons A_{j_1}$. Using the above notations, the creation rate of $c=(c_1,\ldots,c_P)$ reads
\begin{equation}\label{reaction_general}
f(c):=\left(
\begin{array}{c}
 f_1(c)	\\ \vdots \\ f_P(c)
\end{array}
\right)
=\left(
\begin{array}{c|c|c}
\nu_1^1&&\nu_R^1\\
\vdots&\quad\cdots \;\;\quad& \vdots\\
\nu_1^P&	&\nu_R^P
\end{array}
\right)
\left(\hspace{-2mm}
\begin{array}{c}
r_1(c)\\
\vdots\\
r_R(c)
\end{array}\hspace{-2mm}
\right)
=
M
\left(\hspace{-2mm}
\begin{array}{c}
r_1(c)\\
\vdots\\
r_R(c)
\end{array}\hspace{-2mm}
\right).
\end{equation}
Note that the vector field $f$ is quasi-positive. Indeed, we have for all $i\in \{1,\ldots,P\}$,
$$f_i(c)=\sum_{j=1}^R \nu_j^i r_j(c)=\sum_{j:\;\nu_j^i>0}^R \nu_j^i r_j(c)+\sum_{j:\;\nu_j^i<0}^R \nu_j^i r_j(c),$$
and for $c\in \R_+^P$ with $c_i=0$, two cases are possible:
$\nu_j^i>0$ implies $r_j(c)=\nu_j^ik_j^f c_{j_1}c_{j_2}\geq 0$ and
$\nu_j^i<0$ implies $r_j(c)=-\nu_j^ik_j^b c_{j_3}\geq 0$.

Assuming the same diffusion laws as above, the time-evolution of $c=(c_1,\ldots,c_P)$ is now governed by the reaction-diffusion system
\begin{eqnarray}\label{systGeneral}
\left\{
\begin{array}{llll}
\left(
\begin{array}{c}
\partial_t c_1 - \mathrm{div} (d_1(t,x,c)\nabla c_1)\\[0.5ex]
\vdots \\[0.5ex]
\partial_t c_P - \mathrm{div} (d_P(t,x,c)\nabla c_P)\\[0.5ex]
\end{array}
\right)
=\left(
\begin{array}{c}
 f_1(c)	\\ \vdots \\ f_P(c)
\end{array}
\right)\mbox{ on }\; (0,+\infty)\times\Omega,\vspace{2mm}\\
\displaystyle
\partial_\nu c=0\;\mbox{ on }\;(0,+\infty)\times\partial\Omega,\\[1ex]
c(0,\cdot)=c_0\on \Omega.
\end{array}
\right.
\end{eqnarray}
\begin{theo}\label{th2}
Let $p>1$, $s>0$ satisfying $(\ref{intro1})$ and $c_0 \in W_p^s(\Omega,\R_+^P)$. System $(\ref{systGeneral})$ has a unique global nonnegative weak $W^s_p$-solution $c=(c_1,\ldots,c_P):[0,+\infty)\times \Omega \rightarrow \R^P$ provided one of the following conditions is satisfied:
\begin{enumerate}[\quad (i)]
\item $N\leq 3$ and the diffusivities $d_i(t,x,c)$ satisfy $(\ref{hypd_i1})$.
\item $N\leq 5$ and the diffusivities $d_i(c_i)$ satisfy $(\ref{hypd_i2})$.
\end{enumerate}
This solution is actually classical and $(\ref{syst})$ is satisfied in a pointwise sense. Moreover
\begin{equation}\label{Linf:network}
 \forall T>0,\; \exists C=C(\|c_0\|_{L^\infty(\Omega)^P},T )>0 \;\mbox{ such that }\;\|c\|_{L^\infty (Q_T)^P}\leq C.
\end{equation}
If, in addition, $d_i$ and $\partial \Omega$ are smooth, then  $c\in {C}^\infty((0,+\infty)\times \overline\Omega\;;\R_+^P)$.

\end{theo}

As for Theorem~\ref{theomain}, the proof consists in showing that $c$ is uniformly bounded. After deriving a first {\it a priori }estimate from the conservation law $(\ref{conservationlaw})$, or in $L^2(Q_T)$ in the case of diffusivities of the type \eqref{hypd_i2}, we use Lemma \ref{lemEstimates} to improve the regularity of those $c_i$'s whose reaction terms are linearly bounded above. This gives estimates on some quadratic terms and, hence, estimates on some other $c_i$'s. Then we can estimate some further quadratic terms, and so on. Here the atomic conservation law guarantees that we obtain improved estimates for all constituents $c_i$.
Once we have improved the estimates on all the $c_i$'s, we bootstrap this procedure to get estimates in $L^p(Q_T)$ for any $p<+\infty$, and finally in $L^\infty(Q_T)$.

Such a procedure requires that the reactions and the chemical components have been previously sorted. Notice that a permutation of the chemical species corresponds to a permutation of the rows of the stoichiometric matrix $M$, and a permutation of the chemical reactions corresponds to a permutation of its columns.
The concrete way to bring the species and reactions in an appropriate order is based on the following idea: a row in the stoichiometric matrix with only zeros and ones corresponds to a chemical species that is always a product for all of the chemical reactions $A_{j_1}+A_{j_2}\rightarrow A_{j_3}$. If such a species exists, the matrix has a certain block structure.
But, as we assume an atomic mass conservation law, any chemical species whose molar mass is maximal amongst the molar masses of $A_1,\ldots,A_P$ leads to such a row. Indeed, if it would appear as a reactant in $A_{j_1}+A_{j_2}\rightarrow A_{j_3}$, the product $A_{j_3}$ would be heavier - a contradiction.
\begin{lem}\label{lempermutation}
Let $(\ref{conservationlaw})$ be valid. Then, up to a permutation of its rows and columns, the stoichiometric matrix $M$ has the structure
\begin{equation}
\label{merveille}
M={\scriptsize{
\left(\hspace{-1mm}\begin{array}{cccc}
\begin{minipage}[t]{11mm}
\centering{
$\begin{array}[t]{|ccc|}
\hline
 & & \\
 &\hspace{-3mm}\;\;N_1\hspace{-3mm}& \\
 & & \\
& & \\
\hline
1\; & \hspace{-5mm}\;\;\;\ldots\hspace{-5mm} & \hspace{-4mm}1\hspace{-4mm}\\
\hline
\end{array}$}
\end{minipage}
&\hspace{1.5mm}
\begin{minipage}[t]{11mm}
\centering{
$\begin{array}[t]{|ccc|}
\hline
 & & \\
 & & \\
 &\hspace{-3mm}\;\;N_2\hspace{-3mm} & \\
 & & \\
 & & \\
\hline
1 & \hspace{-5mm}\;\;\;\ldots\hspace{-5mm} & \hspace{-4mm}1\hspace{-4mm}\\
\hline
\end{array} $\\[16mm]
\Large{$0$}}
\end{minipage}
 &
 \begin{array}[t]{ccc}
 & & \\
 & & \\
 & & \\
 & & \\
 & & \\
 & & \\
 & & \\
 & & \\
  & & \\
 & \ddots & \\
 & & \\
 & & \\
 & & \vspace{-2mm}\\
 & & \\
  & & \\
 & & \\
\end{array}
&
\begin{minipage}[t]{11mm}
\centering{
$\begin{array}[t]{|ccc|}
\hline
 & & \\
 & & \\
 & & \\
 & & \\
 & & \\
 & & \\
& & \\
 &\hspace{-3mm}\;\;N_k \hspace{-3mm}& \\
 & & \\
& & \\
 & & \\
 & & \\
 & & \\
 & & \\
\hline
1 & \hspace{-5mm}\;\;\;\ldots\hspace{-5mm} & \hspace{-4mm}1\hspace{-4mm}\\
\hline
\end{array} $}
\end{minipage}\hspace{1mm}
\end{array}\;\;\;\right)
}}\;,
\end{equation}
where the submatrices $N_i$ have nonpositive entries.
\end{lem}
{\nd\bf Proof. }
Let $M=(m_{ij})$ and note that, by construction, the columns of $M$ are permutations of the vectors $(-1,-1,+1,0,\ldots,0)$ and $(-2,+1,0,\ldots,0)$. In particular, there is exactly one coefficient equal to $+1$ in each column. Suppose that we have proved the existence of a nonzero row with nonnegative entries. Then, after an appropriate permutation of its rows and columns, $M$ reads
$$
{\scriptsize{
M=\left(\hspace{-1mm}\begin{array}{cc}\vspace{2mm}
\begin{minipage}[t]{11mm}
\centering{
$\begin{array}[t]{|ccc|}
\hline
 & & \\
 & & \\
 &\hspace{-3mm}\;M_1 \hspace{-3mm}& \\
& & \\
& & \\
\hline
0 & \hspace{-5mm}\;\;\;\ldots\hspace{-5mm} & \hspace{-4mm}0\hspace{-4mm}\\
\hline
\end{array}$}
\end{minipage}
&
\begin{minipage}[t]{9mm}
\centering{
$\begin{array}[t]{|ccc|}
\hline
 & & \\
 & & \\
 &\hspace{-3mm}\;N\hspace{-3mm} & \\
 & & \\
 & & \\
\hline
1 & \hspace{-5mm}\;\;\;\ldots\hspace{-5mm} & \hspace{-4mm}1\hspace{-4mm}\\
\hline
\end{array} $}
\end{minipage}
\hspace{1mm}
\end{array}\;\;\;\right)
}}\;,
$$
where $N$ has nonpositive entries and $M_1$ satisfies the same hypothesis as $M$.
By induction, it is then clear that $M$ can be put into the form $(\ref{merveille})$.\saut
Consequently, the proof comes down to finding a nonzero row with nonnegative entries. Let $q\geq1$, $L_{i_1},\ldots,L_{i_q}$ be the rows containing at least one positive entry, and suppose that amongst $L_{i_1},\ldots,L_{i_q}$, every row also has a negative entry. Let $e=(e_1,\ldots,e_P)\in (0,+\infty)^P$ be from $(\ref{conservationlaw})$. By induction, we build a sequence $(u_n)_{n\in \N}$ with values in $\{e_{i_1},\ldots,e_{i_q}\}$ as follows: $u_0=e_{i_1}$;
let $n\geq 0$ and assume that $u_0,\ldots,u_n$ are built such that $u_0<\ldots<u_n$, $u_i\in \{e_{i_1},\ldots,e_{i_q}\}$. By construction, there exist $l\in \{1,\ldots,q\}$ such that $u_n=e_{i_l}$.
The $ i_l^{th}$ row of $M$ has a negative entry by assumption,
so there exists $r\in \{1,\ldots,R\}$ such that $m_{i_lr}\in \{-1,-2\}$. According to $(\ref{conservationlaw})$, the $r^{th}$ column of $M$ satisfies $\langle \nu_r,e\rangle=0$, which reads
$$\exists j\in \{1,\ldots,m\},\exists k\in \{i_1,\ldots,i_q\}:\quad e_{i_l}+e_j=e_k \;\mbox{ if }\;m_{i_lr}=-1,\quad 2e_{i_l}=e_k \;\mbox{ if }\;m_{i_lr}=-2.$$
Then we set $u_{n+1}=e_k$ and, by induction, $(u_n)_{n\in\N}$ is a strictly increasing sequence with values in $\{e_{i_1},\ldots,e_{i_q}\}$: contradiction, so there exists one row amongst $L_{i_1},\ldots,L_{i_q}$ that contains only zeros and ones.\\[-2ex]
$\mbox{ }$ \hfill $\Box$\\
\begin{rem}\label{remestimatesTh2}
According to $(\ref{conservationlaw})$, $P>R$, so the matrix $N_1$ in $(\ref{merveille})$ is nonempty. Let $s\geq 1$ be the number of rows in $N_1$. The point in permuting the rows and columns of $M$ is the following: suppose that $M$ satisfies $(\ref{merveille})$; using the above definition of the reaction terms, there exists $C>0$ depending only on $k_j^f,k_j^b$, such that
\begin{align}
			 1\leq k\leq s\;\Rightarrow \quad &f_k(c)\leq C\sum_{i=1}^P c_i, \label{estTh2a}\\
			 s+1\leq k\leq P\;\Rightarrow \quad &f_k(c)\leq C\left(\sum_{i=1}^P c_i+\sum_{i=1}^{k-1}c_i^2\right).\label{estTh2b}
\end{align}
\end{rem}
{\nd \bf Proof of Theorem \ref{th2}. }As for Theorem \ref{theomain}, the existence of a unique maximal nonnegative weak-$W^s_p$ solution $c=(c_1,\ldots,c_P):[0,T^*)\times \Omega \rightarrow \R^P$ and the regularity results are a consequence of Amann's theory \cite{amann93}. To prove that $T^*=+\infty$, we have to find  bounds in $L^\infty(Q_T)$ for any $T\leq T^*$, $T<+\infty$. Similarly as for Step 1 in the proof of Theorem \ref{theomain}, the first estimates are consequences of the atomic conservation law:
using the no-flux boundary conditions,
\begin{equation*}
\forall t\in (0,T),\quad \sum_{i=1}^P \int_\Omega e_i\, c_i(t)=\sum_{i=1}^P \int_\Omega e_i\, c_{0,i}.
\end{equation*}
Then, using Lemma \ref{lemEstimates} $(i)$, $c$ is bounded in $L^{r_0}(Q_T)^P$ for $r_0\in [1,(N+2)/N)$ if $N\geq 2$, $r_0\in [1,2)$ if $N=1$. For diffusivities $d_i(c_i)$, we write (with $D_i(y)=\int_0^y d_i(s)ds$)
\begin{equation*}
\partial_t\sum_{i=1}^P e_i c_i+\Delta \sum_{i=1}^P e_i D_i(c_i)=0\on Q_T,\;\;
\partial_\nu \sum_{i=1}^P e_iD_i(c_i)=0\on \Sigma_T, \;\;
\sum_{i=1}^P e_ic_i(0,\cdot)=\sum_{i=1}^P e_ic_{0,i}.
\end{equation*}
Then Lemma \ref{lemL2} guarantees that $c$ is bounded in $L^2(Q_T)^P$.\saut
To improve these estimates, using Lemma \ref{lempermutation}, we go down without loss of generality to the case when $M$
has the form given in $(\ref{merveille})$.
Assuming first $N=1$, we know that $c$ is bounded in $L^{r_0}(Q_T)^P$ for $r_0<2$. Using the notations of Remark \ref{remestimatesTh2}, $(\ref{estTh2a})$ and Lemma \ref{lemEstimates} guarantee that $c_1,\ldots,c_s$ are bounded in $L^p(Q_T)$ for any $p<+\infty$. Then, using $(\ref{estTh2b})$, $c_{s+1}$ is bounded in $L^p(Q_T)$ for any $p<+\infty$ and, by induction, for any $k\in \{s+1,\ldots, P\}$, $c_k$ is bounded in $L^p(Q_T)$ for any $p<+\infty$. \saut
Suppose $N\geq 2$ and let $r_0>1$ be such that $c$ is bounded in $L^{r_0}(Q_T)^P$. Using $(\ref{estTh2a})$, $(\ref{estTh2b})$ and Lemma \ref{lemEstimates},
\begin{enumerate}[$\quad$]
\item $c_1,\ldots,c_s$ are bounded in $L^{q_1}(Q_T)$, where $ \frac{1}{r_0}-\frac{2}{N+2}<\frac{1}{q_1}$.
\item $c_1^2,\ldots,c_s^2$ are bounded in $L^{q_2}(Q_T)$, where $ \frac{2}{r_0}-\frac{4}{N+2}<\frac{1}{q_2}$, and $q_2\geq 1$ provided $\frac{2}{r_0}-\frac{4}{N+2}< 1$.
\item $c_{s+1}$ is bounded in $L^{q_3}(Q_T)$, where $\frac{2}{r_0}-\frac{6}{N+2}  <\frac{1}{q_3}$.
\end{enumerate}
Then it is possible to continue improving the estimates for $c_{s+2},\ldots,c_P$ if $q_3\geq q_1$, i.e.\ if
\begin{equation}\label{condTh2improvement}
\frac{2}{r_0}-\frac{6}{N+2}<\frac{1}{r_0}-\frac{2}{N+2}.
\end{equation}
Note that $\frac{2}{r_0}-\frac{4}{N+2}< 1$ is a consequence of $(\ref{condTh2improvement})$. For diffusivities $d_i(t,x,c)$ satisfying $(\ref{hypd_i1})$, $r_0<\frac{N+2}{N}$ and $(\ref{condTh2improvement})$ can be satisfied if and only if $N<4$. For diffusivities $d_i(c_i)$ satisfying $(\ref{hypd_i2})$, $r_0=2$ and  $(\ref{condTh2improvement})$ can be satisfied if and only if $N<6$. Once we have $(\ref{condTh2improvement})$, it is clear that $c_{s+1},\ldots,c_P$ are bounded in $L^{q_1}(Q_T)$ by induction.
Then, similarly as for Theorem \ref{theomain}, we bootstrap this procedure to show that $c$ is bounded in $L^p(Q_T)^P$ for any $p<+\infty$.
\saut
Finally, we use Theorem~\ref{ch5:th} to show that $c_i$ is bounded in $L^\infty(Q_T)$ for all $i$,
whence global existence in Theorem \ref{th2}.\\[-2ex]
$\mbox{ }$ \hfill $\Box$\\[1ex]
\saut
{\nd \bf Example.} For the prototype chain-growth polymerization process, the chemical reaction network reads as
$$A_r+A_1 \underset{k_{r}^b}{\overset{k_r^f}{\rightleftharpoons}} A_{r+1} \;;\; r\in \{1,\ldots,R\},\; k_r^f,k_r^b\geq 0.$$
Typical values for $R$ are large, say about 100 or more. As an example, we write below the equations for $R=4\;$:
\begin{equation*}
\begin{array}{rll}
\left(
\begin{array}{l}
\partial_t c_1-\textrm{div}(d_1(t,x,c)\nabla c_1)\\[0.5ex]
\partial_t c_2-\textrm{div}(d_2(t,x,c)\nabla c_2)\\[0.5ex]
\partial_t c_3-\textrm{div}(d_3(t,x,c)\nabla c_3)\\[0.5ex]
\partial_t c_4-\textrm{div}(d_4(t,x,c)\nabla c_4)\\[0.5ex]
\end{array}
\right)
&=\;\;
\left(
\begin{array}{rrr}
-2&-1&-1\\1&-1&0\\0&1&-1\\0&0&1
\end{array}
\right)
\cdot
\left(
\begin{array}{c}
k_1^f c_1^2-k_1^b c_2\\[0.5ex] k_2^f c_1c_2-k_2^b c_3\\[0.5ex] k_3^f c_1c_3-k_3^b c_4
\end{array}
\right).
\end{array}
\end{equation*}
Note that the stoichiometric matrix is already ``well sorted'' in the sense of Lemma \ref{lempermutation}. Theorem \ref{th2} guarantees the global existence of strong solution for any $R$ in dimension $N=3$ for general diffusivities, and in dimensions $N\leq 5$ for diffusivities $d_i(c_i)$.
Observe that these admissible space dimensions are smaller than for a single reaction
of type
$A_1+A_2
\rightleftharpoons A_3.$

\section{Extensions in the quasi-linear case}
\label{sec:ext}
In the case of diffusion of type $\Delta D_i (c_i)$, some important extensions are
obtained below.
\subsection{Initial data in $\mathbf{L^\infty(\Omega)}$}
For diffusivities of the type $d_i(c_i)$, nonlinear semigroup theory is available to relax the assumption on the initial data in Theorem \ref{th2} to $c_0\in L^\infty(\Omega)_+^P$. A main point is the accretivity of the operator $u\mapsto -\Delta D_i(u)$ with Neumann boundary conditions, which allows to prove uniqueness in the class of bounded mild solutions.
\begin{cor}\label{cor-Linfty-IV}
 Assume $N\leq 3$, $c_0\in L^\infty(\Omega,\R_+^P)$ and let $d_i(c_i)$ satisfy $(\ref{hypd_i1})$. Then system $(\ref{systGeneral})$
 has a unique mild solution $c$ in the class $C([0,+\infty);L^1(\Omega)^P)\cap L^\infty_{loc}([0,+\infty);L^\infty(\Omega)^P)$. Moreover, $c$ is a classical solution on $(0,+\infty)\times \Omega$.
\end{cor}
{\nd\bf Proof.}
For every $\varepsilon >0$, choose $c_0^\varepsilon \in W^2_2 (\Omega)^P$ such that $0\leq c_{0,i}^\varepsilon \leq \|c_{0,i}\|_{L^\infty(\Omega)}$
and $c_0^\varepsilon \to c_0$ a.e.\ on $\Omega$ as $\varepsilon \rightarrow 0$. By Theorem~\ref{th2}, system \eqref{syst} has a unique global classical solution
$c^\varepsilon$, satisfying $(\ref{Linf:network})$. Hence, for all $T>0$, $c^\eps$ is bounded in $L^\infty(Q_T)^P$ independently of $\eps$. Consequently, $f_i(c^\eps)$ is bounded in $L^\infty(Q_T)$. Since $\Omega$ is bounded and $d_i \geq \underline{d}>0$, the semigroup generated by $-\Delta D_i (c_i)$ is compact in $L^1 (\Omega)$. Together, this implies relative compactness of $(c^\varepsilon )_{\varepsilon >0}$ in $C([0,+\infty);L^1 (\Omega)^P)$; cf.\ \cite{Baras}. Then there exists $c$ such that for some sequence $c^k := c^{\varepsilon_k}$, $c^k\to c$ a.e.\ and in $C([0,T ];L^1 (\Omega)^P)$ for any $T>0$. Note in passing that
$c\in L^\infty_{loc}([0,+\infty);L^\infty(\Omega)^P)$. We also have
$f_i (c^k)\to f_i (c)$ a.e.\ and in $L^1 (Q_T)$ for any $T>0$. Consequently, since $c^k$ is also a mild solution, by a standard result on quasi-autonomous evolutions
governed by accretive operators (see, e.g., \cite{BCP} or \cite{ItoKappel}),
$c$ is a mild solution of
\begin{equation}\label{filtration-reaction}
 \partial_t c_i- \Delta D_i (c_i) = f_i (c),\quad \partial_\nu D_i(c_i)_{|\partial \Omega}=0,\quad c_i(0)=c_{0,i}.
\end{equation}
Note also that for locally Lipschitz rate functions $f_i$,
mild solutions of the initial value problem \eqref{filtration-reaction} are unique in the class
$C([0,+\infty);L^1 (\Omega)^P)\cap L^\infty _{loc}([0,+\infty);L^\infty (\Omega)^P)$. Indeed, if $c$ and $\tilde c$ are both mild solutions of  \eqref{filtration-reaction} for the same initial value, then
\begin{equation}\label{gron}
 \sum_{i=1}^P \|c_i (t,\cdot)-\tilde{c}_i (t,\cdot)\|_{L^1 (\Omega)} \leq
\int_0^t \sum_{i=1}^P \|f_i (c(s,\cdot))-f_i (\tilde{c} (s,\cdot))\|_{L^1 (\Omega)} ds.
\end{equation}
Since $f_i$ is Lipschitz continuous on bounded sets and both $c, \tilde c$ are in $L^\infty(Q_T)^P$ for all $T>0$, the Gronwall lemma
guarantees that they coincide.

Let us now prove that $c (t,\cdot) \in W^s_p (\Omega)^P$ for a.e. $t>0$ and for some $s>0$, $p\geq 1$ such that $s>\frac{N}{p}$.
If this holds, we can fix any such $t>0$ and use $c (t,\cdot)$ as the new initial value. By Theorem~\ref{th2},
there exists a classical solution starting at time $t$, and using uniqueness of mild solutions proven above, it coincides with $c$. Hence $c$ is a classical solution on $(t,+\infty)$ for arbitrarily small $t>0$, which ends the proof of Corollary~\ref{cor-Linfty-IV}.\vspace{1mm}
\\
Let $v_i^\varepsilon := D_i (c_i^\varepsilon)$ and note that $v^\varepsilon$ is a solution of
\begin{equation}
\label{d_iDelta}
\partial_t v_i^\varepsilon - d_i(c_i^\varepsilon) \Delta v_i^\varepsilon = d_i(c_i^\varepsilon)\, f_i (c_i^\varepsilon),
\quad \partial_\nu v^\varepsilon_{i|\partial \Omega}=0,\quad v_i(0)=D_i (c_{0,i}).
\end{equation}
Since we do not know whether the coefficient $d_i(c_i^\varepsilon)$ is any better than bounded, we cannot use maximal regularity
in any $L^s$ with large $s$, but have to rely on maximal regularity estimates in $L^2$ as follows.
We consider only a single component $v_i^\varepsilon$.
Let $0<\sigma<\tau<+\infty$.
Multiplication of \eqref{d_iDelta} by $-\Delta v_i^\varepsilon$ and integration over $Q_{\sigma ,\tau}:=[\sigma ,\tau]\times \Omega$
yields
\begin{equation}
\int_\Omega |\nabla v_i^\varepsilon (\tau)|^2 + \underline d\, \int_{Q_{\sigma ,\tau}} |\Delta v_i^\varepsilon |^2 \leq
\int_\Omega |\nabla v_i^\varepsilon (\sigma)|^2 + C_\tau (\tau - \sigma)
\end{equation}
where $C_\tau >0$ is independent of $\varepsilon$. Since $\nabla v_i^\varepsilon=d_i(c_i^\epsilon)\nabla c_i^\epsilon $ and $d_i$ is bounded on bounded sets, there exist $\overline d>0 $ such that for all $\varepsilon>0$,
\begin{equation}
\label{bound1}
\underline d\, \int_{Q_{\sigma ,\tau}} |\Delta v_i^\varepsilon |^2 \leq  {\bar{d}}^{\, 2} \int_\Omega |\nabla c_i^\varepsilon (\sigma)|^2 + C_\tau (\tau - \sigma).
\end{equation}
Next, multiplication of the PDE \eqref{filtration-reaction} for $c_i^\eps$ by $D_i (c_i^\eps)$
and integration over $Q_{\tau}$ yields
\begin{equation}
\int_\Omega \Phi_i (c_i^\varepsilon (\tau)) + \int_{Q_{\tau}} |\nabla D_i (c_i^\varepsilon) |^2
\leq \int_\Omega \Phi_i (c_i^\varepsilon (0)) + C_\tau ,
\end{equation}
again with some constant $ C_\tau >0$, where $\Phi_i (y)=\int_0^y D_i (s) ds$. This implies
\begin{equation}
\int_{Q_\tau} |\nabla c_i^\varepsilon |^2 \leq  K_\tau \quad\mbox{ for some }K_\tau >0.
\end{equation}
In particular, $t\mapsto g^\varepsilon(t):=\int_\Omega |\nabla c_i^\varepsilon (t)|^2 \in L^1 ((0,\tau))_+$.
Hence there exists $\sigma_\varepsilon \in (0,\tau /2)$ such that $g^\varepsilon (\sigma_\varepsilon )\leq 2 |g^\varepsilon|_{L^1 ( (0,\tau/2  ))}/\tau$.
With such $\sigma_\varepsilon$, we obtain by \eqref{bound1} the following estimate on $Q_{\tau/2 ,\tau}$, which is uniform in $\varepsilon >0$:
\[
\underline d\, \int_{Q_{\tau/2 ,\tau}} |\Delta v_i^\varepsilon |^2 \leq \underline d\, \int_{Q_{\sigma_\varepsilon ,\tau}} |\Delta v_i^\varepsilon |^2 \leq 2 \bar d^2\, K_\tau /\tau + \tau C_\tau .
\]
Consequently, for any $i=1,2,3$,
$(\Delta v_i^\varepsilon )_{\varepsilon >0}$ is bounded in $L^2(Q_{\tau/2 ,\tau})$. Then, by \eqref{d_iDelta} and
the uniform $L^\infty$-bounds for the right-hand side on $[0,\tau]$, $(\partial_t v_i^\varepsilon )_{\varepsilon >0}$ is also bounded in $L^2(Q_{\tau/2 ,\tau})$. Hence $(v_i^\varepsilon )_{\varepsilon >0}$ is bounded in $W^1_2((\tau/2,\tau);L^2(\Omega))\cap L^2((\tau/2,\tau);W^2_2(\Omega))$.
By cross-interpolation, $(v_i^\varepsilon )_{\varepsilon >0}$ is bounded in $W^{\eta /2}_2 ( (\tau/2,\tau);W^{2-\eta}_2(\Omega))$ and
the latter is compactly embedded in $L^2 ((\tau/2,\tau)  ;W^{2-\eta}_2(\Omega))$. Hence an appropriate subsequence $(v^{\varepsilon_k}_i)_{k\geq 1}$
converges for a.e. $t\in (\tau/2,\tau)$ in $W^{2-\eta}_2(\Omega)$, and the limit is $D_i(c_i)$.
This shows $D_i (c_i(t))\in W^{2-\eta}_2(\Omega)$ for a.a.\ $t\in [\tau/2,\tau]$, where $\eta >0$ is arbitrarily small. For $p<6$ and $N\leq 3$, we can chose $\eta>0$ small enough so that
$D_i (c_i(t))\in W^1_p(\Omega)$ for a.a.\ $t\in (\tau/2,\tau)$ by Sobolev embedding.
Due to $\nabla c_i = \nabla D_i(c_i)/ d_i (c_i)$ and $\underline d\leq d_i (c_i)$,
this implies $c_i(t) \in W^1_p(\Omega)$ for a.a.\ $t\in [\tau/2,\tau]$.
Since for $p>3$ and $N\leq 3$ we have $1>N/p$ and since $\tau $  can be chosen arbitrarily small, the claim above is proven.\\[-2ex]
$\mbox{ }$\hfill $\Box$
\begin{rem}
Note that in the quasi-linear case, semigroup theory may also be used to prove directly {\it local }well-posedness for system (\ref{systGeneral}) with initial data in $L^\infty(\Omega)_+^P$, see, e.g., \cite{Bo9}.
\end{rem}

\subsection{Rate functions of the type $\mathbf{f(c_1) g(c_2) - h(c_3)}$}
Up to this point, we focused on cases in which the reaction rates read as $k^f c_1 c_2 - k^b c_3$,
corresponding to the chemical transformation
\[
 A + B \rightleftharpoons P.
\]
In concrete applications another prototype case are chemical reactions of type
\[
 A + B \rightleftharpoons P+Q,
\]
leading to rate functions of the form $k^f c_1 c_2 - k^b c_3 c_4$, assuming again elementary reactions, modeled via mass action kinetics.
Now note that the approach from above is not applicable to this case in any space dimension
above one. The reason is the appearance of quadratic terms for both the forward and the backward
reaction path, which destroys the bootstrapping argument.
In fact, $N=2$ is exactly the limiting case for our approach. In the case of constant
diffusivities, it can still be handled, cf.\ \cite{GoudonVasseur}.
For the physical space dimension $N=3$, global existence of strong solutions for this type
of reaction rates is completely open; see \cite{GoudonVasseur, pierre10} for more information.

Due to the quadratic growth, the same problem already appears for reactions of type
\[
 A + B \rightleftharpoons 2P,
\]
also very common in applications.
On the other hand, in Chemical Kinetics there also appear mass action kinetics of fractional orders, obtained via theoretical model simplifications or as empirical rate laws from experimental
measurements. In particular, rate functions regularly contain terms like $c_i^{1/2}$ or
$c_i^{3/2}$ etc.; cf.\ \cite{Espenson, Szabo}.
Therefore, it is also of interest to fix the space dimension to $N=3$ and to investigate,
for instance,
how large the exponent $\gamma$ in a rate function of type $k^f c_1 c_2 - k^b c_3^\gamma$
can be in order to still obtain global existence.

Below, we focus on diffusion operators of type $\Delta D_i (c_i)$. In this case, the $L^2(Q_T)$-bounds are still valid, since they only rely on the cancelation
of the right-hand sides, which makes them a very powerful tool.
For the more general case of $d_i=d_i(t,x,c)$, we refer to Remark~\ref{officialremarkX}.

Since also exponents below 1 are possible, giving more room for the other ones, and since
more elaborate rate functions containing terms like $k_1 c_i + k_2 c_i^{3/2}$ appear (see \cite{Szabo} for this and other examples), we now generalize Theorem~\ref{theomain} to reactions of type
\[
a_1 A_1 + a_2 A_2 \rightleftharpoons a_3 A_3
\]
with constant stoichiometric coefficients $a_i >0$ and associated rate function of the form
\begin{equation}\label{r}
r(c)=r(c_1, c_2, c_3)= f(c_1) g(c_2) - h(c_3),
\end{equation}
where growth conditions on $f,g,h$  have to be imposed so that the uniform estimates $(\ref{Linf:rothe})$ remain valid for the solutions, especially for space dimensions $2$ and $3$.
In the following, we assume
\begin{equation}\label{condition:alpha}
 \left\{
\begin{array}{l}
 f,g,h\in C(\R,\R_+) \;;\; f(0)=g(0)=h(0)=0\;;\vspace{2mm}\\
\ds L:=\underset{s\rightarrow +\infty}
{\limsup} \Big( \frac{f(s)}{s^\alpha}+ \frac{g(s)}{s^\beta}+ \frac{h(s)}{s^\gamma}\Big) <+\infty\; \mbox{ for some }\alpha,\beta,\gamma>0.
\end{array}
\right.
\end{equation}
The following Lemma provides a sufficient condition on $\alpha,\beta,\gamma$ so that the previous $L^\infty$-estimates remain valid on $Q_T$ for any $T>0$.
\begin{lem}\label{lem:alpha}
 Let $c_0\in L^\infty(\Omega)_+^3$, let $T^*\in (0,+\infty]$ and $c:[0,T^*)\times \Omega\rightarrow \R^3$ be a classical solution of (\ref{syst}) with generalized reaction terms
\begin{equation}\label{reac}
F(c):=(-a_1 r(c),-a_2 r(c),a_3 r(c)),
\end{equation}
 where $r$ is defined in \eqref{r} and $f,g,h$ satisfy \eqref{condition:alpha}.
 Assume that $\alpha,\beta,\gamma$ satisfy
\begin{equation}\label{genCondi2}
\begin{array}{c}
\gamma\leq 2 \mbox{ for }N=1\;\;;\;\;
 \gamma \leq 2\mbox{ and }(\alpha+\beta)(\gamma -\frac{4}{N+2})<1+\frac{4}{N+2} \mbox{ for }N\geq 2.
\end{array}
\end{equation}
Then
\begin{equation}\label{Linf:alphacase}
 \forall\; 0<T\leq T^*,\;T<+\infty,\; \exists C=C(\|c_0\|_{L^\infty(\Omega)^3},L,T )>0 \;\mbox{ such that }\;\|c\|_{L^\infty (Q_T)^3}\leq C.
\end{equation}
\end{lem}
{\nd \bf Proof.}
Note that the generalized reaction terms still cancel when considering $F_1/a_1+F_2/a_2+2F_3/a_3$, so Lemma \ref{lemL2} guarantees that $c$ is bounded in $L^2(Q_T)^3$ for any $T>0$. Now the proof follows the lines of what is done for Theorem \ref{theomain}, but Step 2 must be adapted as follows. Let $r_0>1$ such that $c$ is bounded in $L^{r_0}(Q_T)^3$. Using Lemma \ref{lemEstimates}, for $N\geq 2$,
\begin{enumerate}[$\quad$]
\item $h(c_3)$ is bounded in $L^{\frac{r_0}{\gamma}}(Q_T)$, where $\gamma$ satisfies
\begin{equation}\label{cond1bis}
 \gamma \leq r_0.
\end{equation}
\item $c_1,c_2$ are bounded in $L^{q_1}(Q_T)$, where $  \frac{\gamma}{r_0}-\frac{2}{N+2}<\frac{1}{q_1}$;
\item $f(c_1)g(c_2)$ is bounded in $L^{q_2}(Q_T)$, where $  \frac{\gamma(\alpha+\beta)}{r_0}-\frac{2(\alpha+\beta)}{N+2}<\frac{1}{q_2}$. We can choose $q_2\geq1$ provided
\begin{equation}\label{cond2bis}
\frac{\gamma(\alpha+\beta)}{r_0}-\frac{2(\alpha+\beta)}{N+2}< 1;
\end{equation}
\item $c_3$ is bounded in $L^{r_1}(Q_T)$, where $ \frac{\gamma(\alpha+\beta)}{r_0}-\frac{2(\alpha+\beta+1)}{N+2}<\frac{1}{r_1}$.
\end{enumerate}
The initial estimate can be improved if we can choose $r_1>r_0$, i.e.\ if
\begin{equation}\label{cond3bis}
\frac{\gamma(\alpha+\beta)-1}{r_0}<\frac{2(\alpha+\beta+1)}{N+2}.
\end{equation}
If $r_0$ satisfies $(\ref{cond1bis})$, $(\ref{cond2bis})$ and $(\ref{cond3bis})$, the same arguments as in Step 2 in the proof of Theorem $\ref{theomain}$ show that $c$ is bounded in $L^p(Q_T)^3$ for any $p<+\infty$.
Note that this result also holds (with similar computations) for $N=1$, as well as for $r_0$ satisfying $(\ref{cond1bis})$, $(\ref{cond2bis})$ and $(\ref{cond3bis})$ with $N$ replaced by $2$.
Choosing $r_0=2$ in accordance with the initial $L^2(Q_T)$-estimate, $r_0$ satisfies inequalities $(\ref{cond1bis})$, $(\ref{cond2bis})$ and $(\ref{cond3bis})$ if and only if $(\ref{genCondi2})$ holds. Finally, Step 3 carries over to generalized reaction terms without modifications, so Lemma \ref{lem:alpha} holds.
\cqfd
We are now in position to prove the following
\begin{theo}\label{th:alpha}
Let $N=3$ and $c_0\in L^\infty(\Omega)_+^3$. Then system (\ref{syst}) with generalized reaction terms \eqref{reac} where $f,g,h$ satisfy \eqref{condition:alpha} and $\alpha,\beta,\gamma$ satisfy $(\ref{genCondi2})$, has a global mild solution.
\saut
If, in addition, the functions $f,g,h$ are locally Lipschitz continuous, then this solution is unique in the class $C([0,+\infty);L^1(\Omega)^3)\cap L^\infty_{loc}([0,+\infty);L^\infty(\Omega)^3)$. Moreover, $c$ is a classical solution on $(0,+\infty)\times \Omega$.
\end{theo}
{\nd \bf Proof. }
Let $c_0^\varepsilon \in W^2_2 (\Omega)^3$ such that $0\leq c_{0,i}^\varepsilon \leq \|c_{0,i}\|_{L^\infty(\Omega)}$ and  $c_0^\varepsilon \to c_0$ a.e.\ on $\Omega$.  Let $f^\ep,g^\ep,h^\ep\in C^\infty(\R)$ satisfying
(\ref{condition:alpha}) with the same constant $L$
and such that $f^\ep\rightarrow f,g^\ep\rightarrow g,h^\ep\rightarrow h$ as $\ep \rightarrow 0$, uniformly on compact sets. By Amann's results \cite{amann93}, system (\ref{syst}) with initial data $c_0^\ep$ and reaction terms
$F^\ep(c):=(-a_1 r^\ep(c),-a_2 r^\ep(c),a_3 r^\ep(c))$, $ r^\ep(c)=r^\ep(c_1, c_2, c_3)= f^\ep(c_1) g^\ep(c_2) - h^\ep(c_3)$,
has a unique maximal solution $c^\ep:[0,T^*)\times \Omega\rightarrow \R^3$ for some $T^*\in(0,+\infty ]$. By Lemma \ref{lem:alpha}, $ (c^\ep)_{\ep>0} $ is bounded in $L^\infty(Q_T)^3$ for any $T\leq T^*,\;T<+\infty$. In particular, Amann's theory guarantees that $T^*=+\infty$.
Since $c^\ep$ is bounded in $L^\infty(Q_T)^3$ independently of $\ep$, $F^\ep(c^\ep)$ is bounded in $L^\infty(Q_T)^3$ and,
by the compactness result from Baras \cite{Baras}, $c^\varepsilon$ converges, up to extraction of a subsequence,
as $\ep \rightarrow 0$ a.e.\ and in $C([0,T];L^1(\Omega)^3)$ for any $T>0$ to a function $c$. Hence $F^\ep (c^\ep)\rightarrow F(c)$ in $L^1(Q_T)^3$ and therefore $c$ is a global mild solution of (\ref{syst}) by classical nonlinear semigroup theory.

If $f,g,h$ are assumed to be locally Lipschitz continuous, the same Gronwall inequality as (\ref{gron}) in Corollary \ref{cor-Linfty-IV} shows that the mild solution is unique. We may also argue as in the proof of Corollary \ref{cor-Linfty-IV} to show that $c (t,\cdot) \in W^s_p (\Omega)^3$ for a.e.\ $t>0$ and for some $s>0$, $p\geq 1$ such that $s>\frac{N}{p}$. Then H.\ Amann's theory guarantees that there exists a weak $W_s^p$-solution with initial data $c(t)$, which is global by Lemma \ref{lem:alpha} and which is classical on $(t,+\infty)\times \Omega$. Since it is also a mild solution, it coincides with $c$ on $[t,+\infty)\times \Omega$. This is valid for a.e.\ $t>0$, so $c$ is a classical solution on $(0,+\infty)\times \Omega$.
\cqfd

\begin{rem}\label{officialremarkX}
Note that in the case of diffusivites $d_i(t,x,c)$ satisfying $(\ref{hypd_i1})$, we crucially used the linearity in $c_3$ in the reaction terms $c_1c_2-c_3 $ to prove that they are bounded in $L^1(Q_T)$ for any $T>0$, and this was the starting point to derive the first estimate $(\ref{step1})$. For $\gamma>1$, we are no longer able to prove that $f(c_1)g(c_2) -h(c_3) \in L^1(Q_T)$. However, we can do it for $0<\gamma\leq 1$. In the latter case, one can check with similar computations as above that estimate $(\ref{Linf:alphacase})$ still holds for diffusivities $d_i(t,x,c)$ provided
\begin{equation*}
 N=1 \mbox{ or } N\geq 2 \mbox{ and } \gamma <\frac{N+2}{N}\; \mbox{, }\;(\alpha+\beta)(\gamma N-2)<N+2.
\end{equation*}
In the latter case, a similar proof provides the existence of a global mild solution for (\ref{syst}).
\end{rem}


\section{Appendix}\label{S4}
\subsection{Proof of Lemma \ref{lemEstimates}}
{\nd \bf Notations. }Let $\mathcal M =\mathcal M (Q_T,\R)$ be the set of measurable functions on $Q_T$ and for $p\geq 1$, let
\begin{enumerate}[\quad]
\item $L^\infty(0,T;L^p(\Omega))=\{u\in \mathcal M :\underset{{t\in (0,T)}}{\esssup}\|u(t)\|_{L^p(\Omega)}<+\infty\}$, endowed with
$$\|u\|_{L^\infty(0,T;L^p(\Omega))}:=\underset{{t\in (0,T)}}{\esssup}\|u(t)\|_{L^p(\Omega)}\;;$$
\item $L^p(0,T;H^1(\Omega)) = \{u\in \mathcal{ M}:  u\in L^p(0,T;L^2(\Omega))  \mbox{ , }\nabla u \in  L^p(0,T;L^2(\Omega)^N)\}$, endowed with
$$\|u\|_{L^p(0,T;H^1(\Omega))}:= \left( \int_0^T [\|u(t)\|_{L^2(\Omega)}^p +\|\nabla u(t)\|_{L^2(\Omega)^N}^p]dt   \right)^{\frac{1}{p}};$$
\item $V_2(Q_T)=L^\infty(0,T;L^2(\Omega))\cap L^2(0,T; H^1(\Omega))$, endowed with
$$\|u\|_{V_2(Q_T)}:=\left(\|u\|_{{L^\infty(0,T;L^2(\Omega))}}^2 + \|u\|_{L^2(0,T;H^1(\Omega))}^2\right)^{\frac{1}{2}}.$$
\end{enumerate}
\nd To prove Lemma \ref{lemEstimates}, we use the following interpolation result:
\begin{lem}\label{lemmeconsequenceinterpolation}
Let $T>0$, $\Omega$ be a bounded domain of $\R^N$ whose boundary $\partial \Omega$ is at least $C^1$, let $1\leq p<+\infty$ and $u\in L^\infty(0,T;L^p(\Omega))\cap L^2(0,T;H^1(\Omega))$. There exists a constant $C>0$ depending only on $\Omega$, such that
\begin{equation}\label{resultatlemmeinterpol}
\|u\|_{L^q(Q_T)}\leq C\|u\|_{L^\infty(0,T;L^p(\Omega))}^{1-\alpha }\|u\|_{L^2(0,T;H^1(\Omega))}^\alpha ,
\end{equation}
where $\alpha=\frac{2}{q}$ and $q$ satisfies
\begin{equation}\label{leminterpconditions}
 q=2+\frac{2p}{N} \mbox{ for }N\geq 3\;;\quad 2\leq q<2+p \mbox{ for }N=2\;; \quad q=2+p \mbox{ for }N=1.
\end{equation}
\end{lem}
\nd We first recall some classical results: we have the embedding
\begin{equation}\label{1}
 H^1(\Omega)\hookrightarrow L^s(\Omega),
\end{equation}
where $s\geq 1$ satisfies $\frac{1}{s}=\frac{1}{2}-\frac{1}{N} \mbox{ if }N\geq 3 \;;\; s<+\infty \mbox{ if }N=2\;;\; s= +\infty \mbox{ if }N=1$. As a consequence of H\"older's inequality, for $u:\Omega\rightarrow \R$ measurable, $q,r,s\in [1,+\infty]$ and $\alpha \in [0,1]$,
\begin{equation}\label{2}
 \|u\|_{L^q(\Omega)}\leq \|u\|_{L^r(\Omega)}^{1-\alpha}\|u\|_{L^s(\Omega)}^{\alpha},\quad \mbox{where }\;\;\frac{1}{q} =\frac{1-\alpha}{r}+\frac{\alpha}{s}.
\end{equation}
Combining $(\ref{1})$ and $(\ref{2})$, we get the following ``Gagliardo-Nirenberg''-type inequality: there exists $C>0$ depending only on $\Omega$, such that
\begin{equation}\label{3}
  \|u\|_{L^q(\Omega)}\leq C\|u\|_{L^p(\Omega)}^{1-\alpha}\|u\|_{H^1(\Omega)}^{\alpha},
\end{equation}
where $p,q\in [1,+\infty]$, $\alpha \in [0,1]$ and
\begin{equation}\label{4}
 \frac{1}{q}=(1-\alpha)\frac{1}{p} +\alpha(\frac{1}{2}-\frac{1}{N})\mbox{ if }N\geq 3\;;\quad \frac{1-\alpha}{p}<\frac{1}{q} \mbox{ if }N=2\;;\quad \frac{1-\alpha}{p}=\frac{1}{q} \mbox{ if }N=1.
\end{equation}
{\nd\bf Proof of Lemma \ref{lemmeconsequenceinterpolation}.}\\
As $u\in L^\infty(0,T;L^p(\Omega))\cap L^2(0,T;H^1(\Omega))$, we have $u(t)\in L^p(\Omega)\cap H^1(\Omega)$ for a.e.\ $t\in (0,T)$.
Using $(\ref{3})$, we get
\begin{align}
 \int_0^T \|u(t)\|_{L^q(\Omega)}^qdt	 &\leq C^q\int_0^T \|u(t)\|_{L^p(\Omega)}^{q(1-\alpha )}\|u(t)\|_{H^1(\Omega)}^{q\alpha }dt,\nonumber\\
				 &\leq C^q\|u\|_{L^\infty(0,T;L^p(\Omega))}^{q(1-\alpha )}\int_0^T \|u(t)\|_{H^1(\Omega)}^{q\alpha }dt,\label{aab}
\end{align}
where $\alpha $ and $q$ satisfy $(\ref{4})$. Now we choose $q\geq2,\alpha>0$ such that $q\alpha =2$. It is easy to see that conditions $(\ref{4})$ with $q\alpha =2$ are equivalent to conditions $(\ref{leminterpconditions})$. Taking the $(1/q)^{th}$ power in $(\ref{aab})$, we get $(\ref{resultatlemmeinterpol})$.
\cqfd
{\nd \bf Proof of Lemma \ref{lemEstimates}. }\saut
{\it The case $r=1$. }\saut
Integration of $(\ref{eqlemEstimates})$ on $\Omega\times (0,t)$ for $t\in (0,T)$ yields, after integration by parts and using the homogeneous Neumann boundary condition,
\begin{equation}\label{bbbconservationmasse}
 \|u\|_{L^\infty(0,T;L^1(\Omega))}\leq \|f\|_{L^1(Q_T)}+\|u_0\|_{L^1(\Omega)}.
\end{equation}
Let $e=\exp(1)$ and define
$$j:\R_+\rightarrow [0,1),\ y\mapsto 1-\frac{1}{\log{(e+y)}}\quad;\quad J:\R_+\rightarrow \R_+,\;y\mapsto \int_0^y j(s)ds.$$
Multiplication of $(\ref{eqlemEstimates})$ by $j(u)$ and integration by parts on $Q_T$ yields
\begin{align}
\int_\Omega J(u(T))+\int_{Q_T}\frac{d|\nabla u|^2}{(e +u)\log(e+u)^2}&\leq \int_\Omega J(u_0)+\int_{Q_T}fj(u),\nonumber
\end{align}
hence
\begin{align}
\underline d\int_{Q_T}\frac{|\nabla u|^2}{(e+u)\log(e +u)^2}&\leq{\|u_0\|_{L^1(\Omega)}+\|f\|_{L^1(Q_T)}}\label{ap_estimation_grad}.
\end{align}
Let $\beta \in (0,\frac{1}{2})$ and set
$$G:\R_+\rightarrow \R_+,\;y\mapsto \frac{\log(e+y)^2}{(e+y)^{1-2\beta}}\;;\;\|G\|_\infty:=\sup_{y\in\R_+} G(y)<+\infty. $$
Then, for $v=(e+u)^\beta$,
\begin{align}
\int_{Q_T}|\nabla v|^2&=\beta^2\int_{Q_T}\frac{|\nabla u|^2}{(e+u)^{2-2\beta}},\nonumber\\
&=\beta^2\int_{Q_T}\frac{\log(e+u)^2}{(e+u)^{1-2\beta}}\frac{|\nabla u|^2}{(e+u)\log(e+u)^2},\nonumber\\
& \leq \frac{ \|G\|_\infty}{4\underline d}\left({\|u_0\|_{L^1(\Omega)}+\|f\|_{L^1(Q_T)}}\right),\label{aaa}
\end{align}
where we used $(\ref{ap_estimation_grad})$ in the last inequality. According to $(\ref{bbbconservationmasse})$, $v$ is bounded in $L^\infty(0,T;L^{1/\beta}(\Omega))$, so together with $(\ref{aaa})$, $v$ is bounded in $L^\infty(0,T;L^{1/\beta}(\Omega))\cap L^2(0,T;H^1(\Omega))$ and Lemma $\ref{lemmeconsequenceinterpolation}$ guarantees that $v$ is bounded in $L^{r}(Q_T)$, with
$$r=2+\frac{2}{\beta N}\mbox{ for }N\geq 3\;;\quad r<2+\frac{1}{\beta}\mbox{ for }N=2\;;\quad  r=2+\frac{1}{\beta }\mbox{ for }N=1.$$
Then $u$ is bounded in $L^{q}(Q_T)$ with $q=\beta r$, which means
$$q=2\beta+\frac{2}{N}\mbox{ for }N\geq 3\;;\quad q<2\beta+1\mbox{ for }N=2\;;\quad  q=2\beta+1\mbox{ for }N=1.$$
Since $\beta$ can be chosen arbitrarily close to $1/2$, $u$ is bounded in $L^q(Q_T)$, where $q$ satisfies conditions $(i)$ in Lemma \ref{lemEstimates}.
\saut
{\it The case $r>1$. }\saut
Let $p>1$, $t\in (0,T)$. Multiplication of $(\ref{eqlemEstimates})$ by $pu^{p-1}\geq 0$ and integration by parts on $Q_t$ yields
\begin{align}
 \int_{Q_t} \partial_t u^p+ 4(1-\frac{1}{p})\int_{Q_t}d|\nabla (u^{p/2})|^2&\leq p\int_{Q_t}fu^{p-1},\nonumber\\
 \int_\Omega u^p(t)+ 4(1-\frac{1}{p})\int_{Q_t}d|\nabla (u^{p/2})|^2&\leq \int_\Omega u_0^p+p\int_{Q_t}fu^{p-1}\label{l1}.
\end{align}
Here and below, $C$ denotes appropriate constants depending only on $p,\underline d,T$ and $\|u_0\|_{L^\infty(\Omega)}$.
Evidently, $(\ref{l1})$ yields
\begin{equation}\label{appeq1}	
\|u^{p/2}\|_{V_2(Q_{T})}^2 \leq C \left(1+\int_{Q_T}|f|u^{p-1}\right).
\end{equation}
According to Lemma \ref{lemmeconsequenceinterpolation}, we have the continuous embedding $V_2(Q_T)\hookrightarrow L^{s}(Q_T)$, where
\begin{equation}\label{eqboot5}
 s=\frac{2(N+2)}{N}\mbox{ for }N\geq 3\;;\quad s<4\mbox{ for }N=2\;;\quad  s=4\mbox{ for }N=1.
\end{equation}
Assuming $s$ satisfies $(\ref{eqboot5})$, inequality $(\ref{appeq1})$ yields
\begin{equation*}
\exists C>0:\quad \|u^{p/2}\|_{L^{s}(Q_{T})}^2\leq   C \left(1+\int_{Q_T}|f|u^{p-1}\right).
\end{equation*}
Recall that $f\in L^r(Q_T)$, so H\"older's inequality yields
\begin{equation}\label{eqboot2}
 \|u\|_{L^{\frac{ps}{2}}(Q_{T})}^p		\leq C \left(1+ \|f\|_{L^{r}(Q_{T})}\|u\|_{L^{\frac{r(p-1)}{r-1}}(Q_{T})}^{p-1}\right).
\end{equation}
We choose $p>1$ such that
\begin{equation}\label{eqpboostrap2}
1\leq \frac{r(p-1)}{r-1}\leq \frac{ps}{2},
\end{equation}
which is equivalent to
\begin{equation}\label{eqpboostrap22}
1+\frac{s}{2r}-\frac{s}{2}\leq \frac{1}{p}\leq \frac{r}{2r-1}.
\end{equation}
Such a choice is possible if
\begin{equation}\label{eqbootstrap3}
 1+\frac{s}{2r}-\frac{s}{2}<1 \quad \mbox{ and } \quad  1+\frac{s}{2r}-\frac{s}{2}\leq \frac{r}{2r-1}.
\end{equation}
It is easy to check that both inequalities in $(\ref{eqbootstrap3})$ are satisfied for $s\geq 2$, which will be assumed in the following;
note that this is compatible with $(\ref{eqboot5})$. As  $p$ satisfies $(\ref{eqpboostrap2})$, using Young's inequality in $(\ref{eqboot2})$ and $L^{\frac{ps}{2}}(Q_{T})\hookrightarrow L^{\frac{r(p-1)}{r-1}}(Q_{T})$ it follows that
\begin{align}
\exists C>0:\quad\|u\|_{L^{\frac{ps}{2}}(Q_{T})}^{p}&\leq C \left(1+ \|f\|_{L^{r}(Q_{T})}^p+ \frac{1}{2}\|u\|_{L^{\frac{ps}{2}}(Q_{T})}^{p}\right),\label{eqboot3}
\end{align}
and hence $u$ is bounded in $L^{\frac{ps}{2}}(Q_T)$. To get the best estimate, we choose $p$ as large as possible: combining $(\ref{eqboot5})$ with $(\ref{eqpboostrap22})$, we see that the condition on $p$ becomes
\begin{equation}\label{eqboot6}
\frac{N+2}{N}\frac{1}{r}-\frac{2}{N}\leq \frac{1}{p}\mbox{ for }N\geq 3\;;\quad \frac{2}{r}-1<\frac{1}{p}\mbox{ for }N=2\;;\quad  \frac{2}{r}-1\leq \frac{1}{p}\mbox{ for }N=1.
\end{equation}
Since $u$ is bounded in $L^{\frac{ps}{2}}(Q_T)$ with $p$ satisfying $(\ref{eqboot6})$ and $s$ satisfying $(\ref{eqboot5})$, altogether, $u$ is bounded in $L^q(Q_T)$, where $q$ satisfies $(ii)$ in Lemma \ref{lemEstimates}.
\cqfd

\subsection{A priori bounds in $\mathbf{L^\infty}(Q_T)$ for parabolic equations}
\nd In this subsection, we prove that if $c$ satisfies the equation
\begin{equation}\label{ch5:main}
 \left\{
\begin{array}{rcll}
\partial_t c +\divv(-d\nabla c+cu)&=& f	& \on Q_T,\\
-d\partial_\nu c +cu\cdot \nu	&=&0	&\on \Sigma_T,\\
c(0,\cdot)&=&c^0					&\on \Omega,
\end{array}
\right.
\end{equation}
and $f$ is in $L^q(Q_T)$ with $q$ large enough, then $c$ is   bounded in $L^\infty(Q_T)$. This has been shown in \cite{LSU} for the case of Dirichlet boundary conditions (and for general parabolic operators). In the following, we adapt the proof of \cite{LSU} to the case of  Neumann boundary conditions.\vspace{2mm}\\
As before, $\Omega$ is an open, bounded subset of $\R^N$, whose boundary is at least $C^2$. We assume that the data satisfy
\begin{enumerate}[$\quad(i)$]
\item $c^0\in L^\infty(\Omega)_+$.\label{ch5:first}\vspace{1mm}
\item $d:Q_T\rightarrow \R$ is measurable ; $\exists \underline d>0$ such that $\underline d \leq d$.\vspace{1mm}
\item  $|u|^2\in L^\infty(0,T;L^r(\Omega))$ ; $f\in L^{q}(Q_T)\;$ and $\;r,q\geq 1\;$ satisfy
\begin{equation}
 \frac{N}{2r}=1-\theta_1^u\quad;\quad \frac{1}{q}\frac{N+2}{2}=1-\theta_1^f\;,
\end{equation}
where $\theta_1^u,\theta_1^f\in (0,1)$ for $N\geq 2$ and $\theta_1^u,\theta_1^f\in (0,\frac{1}{2})$ for $N=1$.\label{ch5:last}
\end{enumerate}

\begin{theo}\label{ch5:th}
Let $c$ be a classical solution of $(\ref{ch5:main})$ on $Q_T$. Under assumptions $(\ref{ch5:first})-(\ref{ch5:last})$, there exists a constant $M>0$ depending only on the data and $T$, such that
\begin{equation}\nonumber
c(t,x)\leq M \quad\mbox{ for a.e.\ }(t,x)\in Q_T.
\end{equation}
\end{theo}
\nd Let us summarize the notations that will be used in the following:\saut
{\bf Notations. }
Let $c:(0,T)\times \Omega \rightarrow \R$ be a measurable function and $\lambda$ denote the Lebesgue measure, we write
\begin{enumerate}[$\qquad$]
 \item $c_k=\max(0,c-k)$,\quad $k\in \R$.
 \item $Q_{T}(k)=\{ (t,x)\in Q_T: c(t,x)>k \}$.
 \item $A_k(t)=\{x\in \Omega\;:\;c(t,x)>k\}$.
 \item For $q,r\in [1,+\infty]$, the norm on $L^r(0,T;L^q(\Omega))$ is denoted by $\|\cdot\|_{r,q,Q_T}$.
 \item $V_2(Q_T)=L^2(0,T;H^1(\Omega))\cap L^\infty(0,T;L^2(\Omega))$. For $(r,q)$ such that $V_2(Q_T)\inj L^{r}(0,T;L^q(\Omega))$, $\beta>0$ is a constant such that
 $$  \|\cdot\|_{r,q,Q_T}\leq \beta \|\cdot\|_{V_2(Q_T)}. $$
Note that $\beta$ can be chosen independently of $T$ (see \cite{LSU} p. 74).
\end{enumerate}

\nd The subsequent result provides a sufficient condition to deduce uniform bounds on a function from estimates in $V_2(Q_T)$ and $L^{q}(Q_T)$ for finite $q$.
\begin{lem}\label{ch5:lemLSU}
Let $c\in V_2(Q_T)$ and assume that
\begin{equation}\label{ch5:lem:hyp}
\forall k\geq \hat k,\quad \|c_k\|_{V_2(Q_T)}\leq \gamma k \left( \mu_1(k)^{ \frac{1+\theta_1}{r_1}} + \mu_2(k)^{ \frac{1+\theta_2}{r_2}} \right),
\end{equation}
where
$$\hat k,\gamma,\theta_i >0\;;\;\mu_i(k)=\int_0^T \lambda(A_k(t))^{\frac{r_i}{q_i}}\;dt\;;\; \mu_i(k)\in[0,1] \mbox{ for }k\geq \hat k\;;\; i\in \{1,2\},$$
and $(r_i,q_i)$ are chosen such that $V_2(Q_T)\inj L^{r_i}(0,T;L^{q_i}(\Omega)) $. Then there exists $M>0$ depending only on the data, such that for a.e. $(t,x)\in Q_T$,
\begin{equation}\label{ch5:lem:concl}
 c(t,x)\leq M.
\end{equation}
\end{lem}
\nd We will use the following elementary result on numerical sequences:\vspace{2mm}
\begin{enumerate}[$\qquad$]
\item
Let $C,b,\theta >0$ and assume that $(y_n)_{n\in\N}\in \R_+^\N$ satisfies
$$\forall n\in \N,\quad y_{n+1}\leq Cb^n y_n^{\;1+\theta}. $$
Then a straightforward induction on $n$ yields
\begin{equation*}
 \forall n\in \N,\quad y_n \leq C^{\frac{(1+\theta)^n-1}{\theta}} b^{\frac{(1+\theta)^n-1}{\theta^2}-\frac{n}{\theta}}y_0^{(1+\theta)^n}.
\end{equation*}
As a consequence,
\begin{equation}\label{ch5:seq}
\left[b>1 \;\mbox{ and }\; y_0\leq \frac{1}{C^{\frac{1}{\theta}}b^{\frac{1}{\theta^2}}}\right]  \Longrightarrow \;y_n \underset{n\rightarrow +\infty}{\longrightarrow}0.
\end{equation}
\end{enumerate}
{\nd\bf Proof of Lemma \ref{ch5:lemLSU}. } Let $M>\hat k$ , $h_k=M(2-2^{-k})$ for $k\in \N$. It is easy to check that
\begin{equation}\label{ch5:easy}
(h_{k+1}-h_k)\mu_i(h_{k+1})^\frac{1}{r_i}\leq \|c_{h_k}\|_{r_i,q_i,Q_T},\quad i\in \{1,2\}.
\end{equation}
Using $V_2(Q_T)\inj L^{r_i}(0,T;L^{q_i}(\Omega)) $  and $(\ref{ch5:lem:hyp})$,
$$\|c_k\|_{r_i,q_i,Q_T}\leq \beta \|c_k\|_{V_2(Q_T)}\leq \beta\gamma k \left(\mu_1(h_k)^{ \frac{1+\theta_1}{r_1}} + \mu_2(h_k)^{ \frac{1+\theta_2}{r_2}} \right),\quad i\in \{1,2\}. $$
Then
\begin{align}
\mu_i(h_{k+1})^\frac{1}{r_i}&\leq \frac{\|c_{h_k}\|_{r_i,q_i,Q_T}}{h_{k+1}-h_k}\leq \frac{\beta \gamma h_k}{h_{k+1}-h_k}\left( \mu_1(h_k)^{ \frac{1+\theta_1}{r_1}} + \mu_2(h_k)^{ \frac{1+\theta_2}{r_2}}\right )\nonumber\\
&\leq 4\beta \gamma 2^k \left( \mu_1(h_k)^{ \frac{1+\theta_1}{r_1}} + \mu_2(h_k)^{ \frac{1+\theta_2}{r_2}}\right ),\quad i\in \{1,2\}.\label{ch5:l4}
\end{align}
Let $\theta=\min(\theta_1,\theta_2)$. Since $\mu_i(h_k)\in[0,1]$, we have
$$ \mu_1(h_k)^{ \frac{1+\theta_1}{r_1}} + \mu_2(h_k)^{ \frac{1+\theta_2}{r_2}}\leq
\mu_1(h_k)^{ \frac{1+\theta}{r_1}} + \mu_2(h_k)^{ \frac{1+\theta}{r_2}}
\leq C\left(\mu_1(h_k)^{ \frac{1}{r_1}} + \mu_2(h_k)^{ \frac{1}{r_2}}\right)^{1+\theta},$$
where $C>0$ only depends on $\theta$. Going back to $(\ref{ch5:l4})$, we have
$$\mu_1(h_{k+1})^\frac{1}{r_1}+\mu_2(h_{k+1})^\frac{1}{r_2}\leq
8\beta \gamma C2^k\left(\mu_1(h_k)^{ \frac{1}{r_1}} + \mu_2(h_k)^{ \frac{1}{r_2}}\right)^{1+\theta}.
$$
According to $(\ref{ch5:seq})$, the sequence $(\mu_1(h_{k})^\frac{1}{r_1}+\mu_2(h_{k})^\frac{1}{r_2})_{k\in\N}$ converges to $0$ as $k\rightarrow +\infty$ provided its initial value $\mu_1(M)^\frac{1}{r_1}+\mu_2^(M)\frac{1}{r_2}$ is small enough. Similarly as in $(\ref{ch5:easy})$, we have
$$(M-\hat k)\mu_i(M)^{\frac{1}{r_i}}\leq \|c_{\hat k}\|_{r_i,q_i,Q_T},\quad i\in \{1,2\}.$$
Using $(\ref{ch5:lem:hyp})$,
\begin{align*}
 (M-\hat k)  (\mu_1(M)^\frac{1}{r_1}+\mu_2(M)^\frac{1}{r_2})&\leq 2\beta \|c_{\hat k}\|_{V_2(Q_T)}\\
&\leq 2\beta\gamma \hat k \left( \mu_1(h_{\hat k})^\frac{1+\theta_1}{r_1}+\mu_2(h_{\hat k})^\frac{1+\theta_2}{r_2}\right)\\
&\leq 4\beta \gamma \hat k \mbox{ \quad (since $\mu_i(h_{\hat k}) \in [0,1]$). }
\end{align*}
We deduce that $\mu_1(M)^\frac{1}{r_1}+\mu_2(M)^\frac{1}{r_2}$ can be chosen arbitrarily small provided $M$ is large enough, and then
$$\mu_1(2M)^\frac{1}{r_1}+\mu_2(2M)^\frac{1}{r_2}\leq \mu_1(h_{k})^\frac{1}{r_1}+\mu_2(h_{k})^\frac{1}{r_2}\underset{k\rightarrow +\infty}{\longrightarrow}0,$$
whence $c(t,x)\leq 2M$ for a.e. $(t,x)\in Q_T$.
\cqfd

{\nd \bf Proof of Theorem \ref{ch5:th}. }
Let $k\geq \|c^0\|_{L^\infty(\Omega)}$. We multiply equation $(\ref{ch5:main})$ by $ c_k$, integrate on $Q_{t_1}$ for $t_1\in (0,T)$ and integrate by parts to get, using the homogeneous Neumann boundary conditions,
\begin{align*}
 \int_{Q_{t_1}(k)}\frac{1}{2} \partial_t (c_k^2 )+\int_{Q_{t_1}(k)}d|\nabla c_k|^2 &=\int_{Q_{t_1}(k)}c\;u\cdot \nabla c_k +fc_k\;,\\
 \frac{1}{2}\int_\Omega c_k^2(t_1) +\underline d\int_{Q_{t_1}(k)}|\nabla c_k|^2 &\leq \int_{Q_{t_1}(k)}|c|\;|u|\;|\nabla c_k| +|f|c_k.
\end{align*}
Using Young's inequality to absorb the term $\nabla c_k$ in the left-hand side, there exists $\alpha=\alpha(\underline d)>0$ such that
$$
\alpha  \left[\int_\Omega c_k^2(t_1) +\int_{Q_{t_1}(k)}|\nabla c_k|^2\right] \leq \int_{Q_{t_1}(k)} |u|^2|c|^2 +|f|c_k,
$$
and consequently
\begin{equation*}
\alpha  \|c_k\| _{V_2(Q_{t_1} )}^2  \leq \int_{Q_{t_1}(k)} |u|^2|c|^2 +|f|c_k.
\end{equation*}
From now on, we impose $k\geq 1$, so that
\begin{equation}\label{ch5:l1}
\alpha  \|c_k\| _{V_2(Q_{t_1} )}^2  \leq   \int_{Q_{t_1}(k)} |u|^2|c|^2 +|f|c_k \leq 2\int_{Q_{t_1}(k)} (|u|^2+|f|) (c_k^2+k^2).
\end{equation}
We now estimate the right-hand side as follows:
\begin{align*}
 \int_{Q_{t_1}(k)}|u|^2 (c_k^2+k^2) &\leq \||u|^2  \|_{{\infty,r,Q_{t_1}(k)}} \|c_k^2+k^2\|_{{1,\frac{r}{r-1},Q_{t_1}(k)  }  }\\
			      &\leq \||u|^2  \|_{{\infty,r,Q_{t_1}(k)}}\big(  \|c_k\|^2_{{2,\frac{2r}{r-1},Q_{t_1}(k)  } } +k^2 \|1\|_{{1,\frac{r}{r-1},Q_{t_1}(k)  } } \big).
\end{align*}
Using  H\"older's inequality,
$$
\|c_k\|_{2, \overline r,Q_{t_1}(k)  }\leq \|c_k\|_{2(1+\theta^u)  ,\hat r,Q_{t_1}(k)}\; \mu_u(k)^{\frac{\theta^u }{2 (1+\theta^u)} },
$$
where
$$\mu_u(k)=\int_0^{t_1} \lambda(A_k(t))^{ \frac{r-1}{r}  } dt\;;\;  \overline r=\frac{2r}{r-1} \;;\; \hat r = \overline r (1+\theta^u)\;;\; \theta^u= \frac{2\theta_1^u}{N}.$$
It is easy to check that
$$\frac{1}{2 }+ \frac{N}{2\overline{r}}=\frac{N}{4} +\frac{\theta_1^u}{2}\quad;\quad \frac{1}{2(1+\theta^u)}+\frac{N}{2 \hat r}=\frac{N}{4},$$
and therefore we have the embedding $V_2(Q_{t_1})\inj L^{2(1+\theta^u)} (0,t_1;L^{\hat r}(\Omega))$ (see e.g. \cite{LSU} p.74). As a consequence, there exists $\beta>0$ (independent of $t_1$), such that
\begin{equation}\label{ch5:l1a}
 \|c_k\|^2_{{2,\overline{r},Q_{t_1}(k)  } }\leq \beta^2 \|c_k\|_{V_2(Q_{t_1})}^2\mu_u(k)^{\frac{\theta^u }{1+\theta^u}}.
\end{equation}
For the second term, we have
\begin{equation}\label{ch5:l1b}
k^2 \|1\|_{{1,\frac{r}{r-1},Q_{t_1}(k)  } }=k^2 \left( \int_0^{t_1} \lambda(A_k(t))^{\frac{r-1}{r}} dt  \right)=k^2 \mu_u(k).
\end{equation}
Similarly,
\begin{align*}
 \int_{Q_{t_1}(k)}|f| (c_k^2+k^2) &\leq \|f  \|_{{q,Q_{t_1}(k)}} \big( \|c_k^2+k^2\|_{{\frac{q}{q-1},Q_{t_1}(k)  }  } \big)\\
			&\leq \|f  \|_{{q,Q_{t_1}(k)}} \big (\|c_k\|^2_{{\frac{2q}{q-1},Q_{t_1}(k)  }  } +k^2 \|1\|_{_{{\frac{q}{q-1},Q_{t_1}(k)  }  }}
			\big).
\end{align*}
Then using H\"older's inequality,
$$
 \|c_k\|_{{\overline{q},Q_{t_1}(k)}}\leq \|c_k\|_{{\hat q,Q_{t_1}(k)}} \mu_f(k)^{\frac{1}{\overline q}-\frac{1}{\hat q}},
$$
where
$$\mu_f(k)=\int_0^{t_1}  \lambda(A_k(t))dt\;;\; \overline{q}=\frac{2q}{q-1}\;;\; \hat q=\overline q(1+\theta^f)\;;\;\theta^f=\frac{2\theta_1^f}{N}.$$
One can check that $\frac{1}{\overline q}+\frac{N}{2\overline q}=\frac{N}{4}+\frac{\theta_1^f}{2}$, $\frac{1}{\hat q}+\frac{N}{2\hat q}=\frac{N}{4}$, so $V_2(Q_{t_1})\inj L^{\hat q}(Q_{t_1} )$ and therefore
\begin{equation}\label{ch5:l1c}
\|c_k\|^2_{{\frac{2q}{q-1},Q_{t_1}(k)}}\leq \beta^2 \|c_k\|_{V_2(Q_{t_1})}^2\mu_f(k)^{\frac{2\theta^f}{\hat q}}.
\end{equation}
The last term is
\begin{equation}\label{ch5:l1d}
 k^2 \|1\|_{_{{\frac{q}{q-1},Q_{t_1}(k)  }  }}=k^2 \mu_f(k)^{\frac{2(1+\theta^f)}{\hat q}   }.
\end{equation}
Going back to $(\ref{ch5:l1})$ and using $(\ref{ch5:l1a})-(\ref{ch5:l1d})$, there exists $C>0$ depending only on $\beta$, \linebreak $\||u|^2\|_{{\infty,r,Q_{T}(k)}}$ and $\|f\|_{L^q(Q_T)}$ (but not on $t_1$), such that for all $k\geq \max( \|c^0\|_{L^\infty(\Omega)},1)$,
\begin{equation}\label{ch5:l2}
 \alpha  \|c_k\|_{V_2(Q_{t_1} )}^2  \leq C \left[\|c_k\|_{V_2(Q_{t_1})}^2 \big(\mu_u(k)^{\frac{\theta^u}{1+\theta^u}}+ \mu_f(k)^{\frac{2\theta^f}{\hat q}}  \big)+k^2 \big(\mu_u(k)+\mu_f^{\frac{2(1+\theta^f)}{\hat q}}(k)\big)\right].
\end{equation}
We now choose $t_1\in (0,T)$ small enough so that
\begin{equation*}
 \ds C  \big(\mu_u(k)^{\frac{\theta^u}{1+\theta^u}}+ \mu_f(k)^{\frac{2\theta^f}{\hat q}}  \big)\leq \frac{\alpha}{2}
\quad ;\quad  t_1 \lambda (\Omega)^{\frac{r-1}{r}}\leq 1\quad;\quad t_1\lambda(\Omega)\leq 1.
\end{equation*}
This is the case provided
\begin{equation}\label{ch5:l3}
 C \big(t_1^{\frac{\theta^u}{1+\theta^u}} \lambda(\Omega)^{\frac{2\theta^u}{\hat r} }+
t_1^{\frac{2\theta^f}{\hat q}} \lambda(\Omega)^{\frac{2\theta^f}{\hat q} }
\big)\leq \frac{\alpha}{2}\quad ;\quad  t_1 \lambda (\Omega)^{\frac{r-1}{r}}\leq 1\quad;\quad t_1\lambda(\Omega)\leq 1.
\end{equation}
For $t_1$ satisfying $(\ref{ch5:l3})$, inequality $(\ref{ch5:l2})$ yields, if $\hat k=\max( \|c^0\|_{L^\infty(\Omega)},1)$,
\begin{equation}
\forall k\geq \hat k,\quad  \frac{\alpha }{2} \|c_k\|_{V_2(Q_{t_1} )}^2  \leq k^2C \big(\mu_u(k)^{\frac{2(1+\theta^u)}{2(1+\theta^u)}} +\mu_f(k)^{\frac{2(1+\theta^f)}{\hat q}}\big).
\end{equation}
Moreover, for all $k\geq \hat k$, $\;\mu_u(k), \mu_f(k)\in [0,1]$, so we can apply Lemma \ref{ch5:lemLSU} and $c$ is bounded on $Q_{t_1}$. Remark that $t_1$ does not depend on $\hat k$. Then we may subdivide $Q_T=(0,T)\times\Omega$ in a finite sequence of cylinders $(t_i,t_{i+1})\times\Omega $, $i=1\ldots,P$, whose altitudes $(t_{i+1}-t_i)$ are subject to the requirement $(\ref{ch5:l3})$. Applying the above result on each cylinder, we get that $c$ is bounded on $Q_T$, which ends the proof of Theorem \ref{ch5:th}.

\bibliography{biblio}{}

\begin{thebibliography}{dGM84}

\bibitem[Ama89]{amann89}
H.~Amann.
\newblock Dynamic theory of quasilinear parabolic systems. {III}. {G}lobal
  existence.
\newblock {\em Math. Z.}, 202(2):219--250, 1989.

\bibitem[Ama93]{amann93}
H.~Amann.
\newblock Nonhomogeneous linear and quasilinear elliptic and parabolic boundary
  value problems.
\newblock In {\em Function spaces, differential operators and nonlinear
  analysis ({F}riedrichroda, 1992)}, volume 133 of {\em Teubner-Texte Math.},
  pages 9--126. Teubner, Stuttgart, 1993.

\bibitem[Bar78]{Baras}
P.~Baras.
\newblock {Compacit\'e de l'op\'erateur $f\to u$ solution d'une \'equation non
  lin\'eaire $(du/dt)+Au\ni f$}.
\newblock {\em C. R. Acad. Sci Paris.}, 268:1113--1116, 1978.

\bibitem[BB99]{BaldygaBourne}
J.~Baldyga and J.R. Bourne.
\newblock {\em Turbulent mixing and chemical reactions}.
\newblock John Wiley \& Sons, New York, 1999.

\bibitem[BCP]{BCP}
Ph. B\'enilan, M.G. Crandall, and A.~Pazy.
\newblock {\em Evolution Equations in Banach Spaces, preprint monograph}.

\bibitem[BD]{BotheDreyer}
D.~Bothe and W.~Dreyer.
\newblock Continuum thermodynamics of chemically reacting fluid mixtures.
\newblock {\em arXiv:1401.5991, preprint 2014}.

\bibitem[Bot96]{Bo9}
D.~Bothe.
\newblock Flow invariance for perturbed nonlinear evolution equations.
\newblock {\em Abstract and Applied Analysis}, 1:379--395, 1996.

\bibitem[Bot99]{BoHabil}
D.~Bothe.
\newblock {\em Nonlinear Evolutions in Banach Spaces -- Existence and
  Qualitative Theory with Applications to Reaction-Diffusion Systems}.
\newblock Habilitation thesis, University of Paderborn, 1999.

\bibitem[Bot11]{BoMS}
D.~Bothe.
\newblock On the {Maxwell-Stefan} equations to multicomponent diffusion.
\newblock In {\em Progress in Nonlinear Differential Equations and their
  Applications}, volume~60, pages 81--93. Springer, Basel, 2011.

\bibitem[BP10]{BP-intermediate}
D.~Bothe and M.~Pierre.
\newblock Quasi-steady-state approximation for a reaction-diffusion system with
  fast intermediate.
\newblock {\em J. Math. Anal. Appl.}, 368:120--132, 2010.

\bibitem[BPR12]{BPR12}
D.~Bothe, M.~Pierre, and G.~Rolland.
\newblock Cross-diffusion limit for a reaction-diffusion system with fast
  reversible reaction.
\newblock {\em Comm.\ in Part.\ Diff.\ Eqs.}, 37:1940--1966, 2012.

\bibitem[Cus97]{cussler}
E.L. Cussler.
\newblock {\em Diffusion: {M}ass {T}ransfer in {F}luid {S}ystems}.
\newblock Cambridge University Press, 2 edition, 1997.

\bibitem[DF06]{DF06}
L.~Desvillettes and K.~Fellner.
\newblock Exponential decay toward equilibrium via entropy methods for
  reaction-diffusion equations.
\newblock {\em Journal of Mathematical Analysis and Applications},
  319:157--176, 2006.

\bibitem[dGM84]{deGrootMazur-book}
S.R. de~Groot and P.~Mazur.
\newblock {\em Non-Equilibrium Thermodynamics}.
\newblock Dover Publications, 1984.

\bibitem[DHP03]{denk}
R.~Denk, M.~Hieber, and J.~Pr{\"u}ss.
\newblock {$\RR$}-{B}oundedness, {F}ourier {M}ultipliers and {P}roblems of
  {E}lliptic and {P}arabolic type.
\newblock {\em Mem. Amer. Math. Soc.}, 166(788):1--114, 2003.

\bibitem[Esp95]{Espenson}
J.H. Espenson.
\newblock {\em Chemical Kinetics and Reaction Mechanisms}.
\newblock McGraw-Hill, 2 edition, 1995.

\bibitem[{\'E}T89]{erdi}
P.~{\'E}rdi and J.~T{\'o}th.
\newblock {\em Mathematical models of chemical reactions}.
\newblock Nonlinear Science: Theory and Applications. Princeton University
  Press, Princeton, NJ, 1989.
\newblock Theory and applications of deterministic and stochastic models.

\bibitem[Fen91]{feng91}
W.~Feng.
\newblock Coupled system of reaction-diffusion equations and applications in
  carrier facilitated diffusion.
\newblock {\em Nonlinear Anal.}, 17(3):285--311, 1991.

\bibitem[Gio99]{Giovan}
V.~Giovangigli.
\newblock {\em Multicomponent Flow Modeling}.
\newblock Birkhäuser, Boston, 1999.

\bibitem[GV10]{GoudonVasseur}
T.~Goudon and A.~Vasseur.
\newblock Regularity analysis for systems of reaction-diffusion equations.
\newblock {\em Ann. Sci. Ec. Norm. Sup.}, 43:117--142, 2010.

\bibitem[IK02]{ItoKappel}
K.~Ito and F.~Kappel.
\newblock {\em {Evolution Equations and Approximation}}.
\newblock World Scientific, 2002.

\bibitem[JS13]{Juengel}
A.~Jüngel and I.~V. Stelzer.
\newblock Existence analysis of {Maxwell-Stefan} systems for multicomponent
  mixtures.
\newblock {\em SIAM J. Math. Anal.}, 2013.

\bibitem[KT93]{TaylorKrishna}
R.~Krishna and R.~Taylor.
\newblock {\em Multicomponent mass transfer}.
\newblock Wiley (New York), 1993.

\bibitem[LSO96]{transport-in-pm}
P.C. Lichtner, C.~I. Steefel, and E.H. Oelkers.
\newblock {\em Reactive Transport in Porous Media}, volume~34 of {\em Reviews
  in Mineralogy}.
\newblock Mineralogical Society of America, 1996.

\bibitem[LSU68]{LSU}
O.A. Ladyzenskaja, V.A. Solonnikov, and N.N. Ural'seca.
\newblock {\em Linear and quasilinear equations of parabolic type}, volume~23
  of {\em Translation of {M}athematical {M}onographs}.
\newblock Amer. Math. Soc, 1968.

\bibitem[Mur89]{Murray}
J.D. Murray.
\newblock {\em Mathematical Biology}.
\newblock Springer, Berlin, 1989.

\bibitem[MW04]{MorganWaggonner}
J.~Morgan and S.~Waggonner.
\newblock Global existence for a class of quasilinear reaction-diffusion
  systems.
\newblock {\em Commun. Appl. Anal.}, 8:153--166, 2004.

\bibitem[Pie10]{pierre10}
M.~Pierre.
\newblock Global existence in reaction-diffusion systems with control of mass:
  a survey.
\newblock {\em Milan J. Math.}, 78:417--455, 2010.

\bibitem[Rot84]{rothe}
F.~Rothe.
\newblock {\em Global solutions of reaction-diffusion systems}, volume 1072 of
  {\em Lecture Notes in Mathematics}.
\newblock Springer-Verlag, Berlin, 1984.

\bibitem[Sza64]{Szabo}
Z.G. Szabo.
\newblock {\em {Advances in Kinetics of Homogeneous Gas Reactions}}.
\newblock Methuen \& Co Ltd., London, 1964.

\bibitem[Vaz07]{Vazquez}
J.L. Vazquez.
\newblock {\em The Porous Medium Equation - Mathematical Theory}.
\newblock Clarendon-Press, Oxford, 2007.

\end{thebibliography}
\bibliographystyle{alpha}
\end{document}